\pdfoutput=1
\documentclass[11pt,letterpaper]{article}
\usepackage[T1]{fontenc}
\usepackage[utf8]{inputenc}
\usepackage[margin=1in]{geometry}
\usepackage{amsmath,amssymb,amsthm,txfonts,bm}
\usepackage{graphicx,booktabs,array}
\usepackage{enumitem}
\usepackage[section]{placeins}
\usepackage[round,authoryear]{natbib}
\bibpunct{(}{)}{;}{a}{}{,}
\usepackage{microtype}
\usepackage[hidelinks,hypertexnames=false]{hyperref}
\hypersetup{
  pdftitle={Neyman Orthogonalized Specification Tests for Conditional Tail Index Models},
  pdfauthor={Liuhua Peng, Shuo Li, and Xiaojun Song},
  pdfsubject={Neyman orthogonalization; conditional tail index specification; second-order regular variation}}
\providecommand{\doi}[1]{\href{https://doi.org/#1}{doi:\nolinkurl{#1}}}
\allowdisplaybreaks[2]

\newcommand{\X}{\boldsymbol X}
\newcommand{\x}{\boldsymbol x}
\newcommand{\tht}{\boldsymbol\theta}
\newcommand{\vv}{\boldsymbol v}
\newcommand{\dd}{\boldsymbol d}

\newcommand{\trans}{^{\mathsf T}}
\newcommand{\ind}{\mathbf 1}
\newcommand{\E}{\mathbb E}
\newcommand{\Prob}{\mathbb P}
\newcommand{\R}{\mathbb R}
\newcommand{\op}{o_{p}}
\newcommand{\Op}{O_{p}}

\newcommand{\pn}{\mathbb P_N}

\newtheorem{theorem}{Theorem}[section]
\newtheorem{proposition}[theorem]{Proposition}
\newtheorem{corollary}[theorem]{Corollary}

\theoremstyle{definition}
\newtheorem{assumption}{Assumption}

\newtheorem{example}[theorem]{Example}
\theoremstyle{remark}
\newtheorem{remark}[theorem]{Remark}
\numberwithin{equation}{section}

\title{Neyman Orthogonalized Specification Tests for Conditional Tail Index Models}
\author{Liuhua Peng\thanks{\raggedright School of Mathematics and Statistics, The University of Melbourne,
  Victoria 3010, Australia. \mbox{E-mail: \href{mailto:liuhua.peng@unimelb.edu.au}{\nolinkurl{liuhua.peng@unimelb.edu.au}}}.}
  \and Shuo Li\thanks{\raggedright School of Statistics, Tianjin University of Finance and Economics,
  Tianjin 300222, China. \mbox{E-mail: \href{mailto:shuoli@tjufe.edu.cn}{\nolinkurl{shuoli@tjufe.edu.cn}}}.}
  \and Xiaojun Song\thanks{\raggedright Department of Business Statistics and Econometrics,
  Guanghua School of Management and Center for Statistical Science, Peking University,
  Beijing, 100871, China. \mbox{E-mail: \href{mailto:sxj@gsm.pku.edu.cn}{\nolinkurl{sxj@gsm.pku.edu.cn}}}.}}
\date{}

\begin{document}

\maketitle
\begin{abstract}
    We develop specification tests for parametric conditional tail index models of heavy-tailed distributions.
    The tests are based on a projected empirical process that removes the first-order effect of estimating the null parameter.
    Within the conditional Hall class, we identify conditions under which projection also eliminates the leading Pareto approximation bias in the test process, permitting lower thresholds while preserving asymptotically correct size. 
    We establish the validity of the tests via multiplier bootstrap and characterize power against fixed and local alternatives.
    For log-linear models, bootstrap validity extends to settings in which covariates of observations above a common threshold concentrate and the smallest eigenvalue of the information matrix tends to zero. 
    We further propose a novel threshold selection rule based on an estimate of the Pareto approximation bias after projection. 
    The rule uses the estimated bias to determine whether a reference threshold can be lowered to retain more observations, with potential gains in power. 
    The finite sample performance of the tests is examined through simulations and a real data application.
\end{abstract}

\noindent\textit{Keywords:} Conditional Hall class; Empirical process; Multiplier bootstrap; Pareto approximation; Threshold selection.

\section{Introduction}
\label{sec:introduction}

The tail index is a fundamental quantity in extreme value theory \citep{dehaan2006}. 
For heavy-tailed distributions, it governs the rate at which tail probabilities decay and provides the basis for extrapolating extreme quantiles beyond the observed range \citep{weissman1978}. 
In applications, this rate may depend on observable characteristics. 
For example, insurance claim sizes may have heavier tails in some lines of business than in others, while the tail behavior of financial losses may vary with market conditions. 
The conditional tail index characterizes differences in tail decay across covariate values that need not be apparent from the conditional mean.

Nonparametric methods estimate the conditional tail index without specifying a functional form for its dependence on the covariates \citep{girard2022}. 
Semiparametric alternatives include partially linear models \citep{li2022}, varying-coefficient models \citep{momoki2024}, and single-index models \citep{yoshida2025}. 
Parametric models are widely used in applications, often with the logarithm of the conditional tail index specified as a linear function of the covariates \citep{wangtsai2009,nicolau2023}. 
Such models are particularly attractive when observations above a high threshold are scarce, since observations at different covariate values contribute to the estimation of a common model parameter. 
A misspecified model, however, may imply an incorrect rate of tail decay at some covariate values. 
In this paper, we develop tests of the null hypothesis that the conditional tail index belongs to a specified parametric family.

A natural approach is to compare the conditional tail index fitted under the null parametric model with a nonparametric kernel estimator. 
Both estimates rely on the relatively few observations above a high threshold, and kernel smoothing further reduces the effective sample size. 
Under the null hypothesis, the difference between the two estimators therefore converges at the slower nonparametric rate. 
A test based on these differences can have low power against local departures from the null hypothesis, particularly when the dimension of the covariate is large.

To avoid nonparametric estimation of the conditional tail index, we construct an empirical process based on the log excess and the null model. 
Let $Y$ denote the outcome, $\boldsymbol{X}$ the covariate vector, and $Z=\log(Y/u)$ the log excess when $Y$ exceeds a high threshold $u$. 
Under the Pareto approximation, $Z$ conditional on $Y>u$ and $\boldsymbol{X}=\boldsymbol{x}$ is exponentially distributed with mean equal to the conditional tail index at $\boldsymbol{x}$. 
At each log excess level $t$, we form a residual by subtracting the approximated conditional probability of $Z>t$ under the null parametric model from the indicator $\mathbf{1}\{Z>t\}$. 
Following the integrated conditional moment approach of \citet{bierens1990} and \citet{stute1997}, we sum these residuals over covariate regions $\{\boldsymbol{X}\leq\boldsymbol{x}\}$. 
The resulting process provides a basis for testing the conditional tail index specification without local smoothing.

Estimating the model parameter affects the empirical process through the approximated conditional probabilities of $Z>t$. 
Even under an exact Pareto model, the expansion of the process generally contains a term linear in the parameter estimation error. 
The error in the Pareto approximation can also bias the parameter estimate and shift the mean of the test process even when the tail index model is correctly specified. 
To remove the effect of parameter estimation, we construct a Neyman-orthogonalized process \citep{neyman1959,chernozhukov2018} through projection. 
We project the indicators of covariate regions onto the derivatives of the approximated conditional probabilities of $Z>t$ with respect to the model parameter. 
The projection residuals weight the differences between $\mathbf{1}\{Z>t\}$ and the approximated conditional probabilities and form the projected test process. 
Orthogonality eliminates the first-order contribution of parameter estimation to the test process. 
We form Cram\'er--von Mises (CvM) and Kolmogorov--Smirnov (KS) statistics from the projected test process and obtain critical values by a multiplier bootstrap.

Implementing the proposed tests requires selecting a threshold above which observations are used. 
The choice of threshold involves a trade-off between the accuracy of the Pareto approximation and the number of observations available for testing \citep{bader2018,schneider2021}. 
Raising the threshold generally improves the Pareto approximation but leaves fewer observations with which to detect departures from the null model. 
Lowering the threshold retains more observations but may increase the Pareto approximation bias in the test process.
If this bias is not asymptotically negligible, the null rejection probabilities of the proposed tests can exceed the nominal level.

Interestingly, we find that in some conditional tail models the projection used to remove the first-order effect of parameter estimation also eliminates the leading Pareto approximation bias in the test process. 
At each log excess level $t$, this occurs when the leading approximation error in the conditional survival probability lies in the span of its derivatives with respect to the model parameter.
Examples include conditional tail index models constructed from Burr, Student-$t$, and Fr\'echet distributions. 
The remaining bias is then of smaller order, allowing the multiplier bootstrap to approximate the null distributions of the CvM and KS statistics consistently at lower thresholds.

Motivated by this finding, we propose a novel threshold selection rule for the proposed tests.
At each candidate threshold, we estimate the Pareto approximation bias after projection by comparing the conditional means of the projected test process under an estimated conditional distribution of the log excess and an exponential distribution whose mean is an estimate of the conditional tail index computed without imposing the null model. 
For each test, we select the lowest candidate threshold satisfying a prescribed tolerance, retaining a reference threshold if no candidate qualifies.

This paper makes several contributions. 
First, we develop Neyman-orthogonalized specification tests for parametric conditional tail index models, allowing the index to vary with the covariates under the null. 
Related tests concern tail-index homogeneity \citep{yanghomogeneity2026} or varying-coefficient models \citep{momoki2024}. 
Empirical process tests in extreme value theory have been developed for domains of attraction \citep{drees2006tail} and tail copulas \citep{can2015tail}, but to our knowledge this is the first use of a projected empirical process to test a conditional tail index specification. 
We establish asymptotically correct size for the proposed tests and characterize their asymptotic power against fixed and local alternatives. 
As a common threshold increases, the covariates of observations
with $Y_i>u$ can concentrate near values where the tail index
is largest \citep{yangrank2026}. 
For log-linear models, we establish bootstrap validity under conditions that allow the smallest eigenvalue of the information matrix to tend to zero as a result of this concentration. 

Second, we identify an additional role for Neyman orthogonalization.
Projection removes the first-order effect of parameter estimation
in specification tests for quantile and propensity score models
\citep{escancianogoh2014,santanna2019,escanciano2026double}.
In conditional tail index tests, the Pareto approximation can shift the mean of the projected test process. 
For the conditional Hall class \citep{hall1982}, the leading term in the mean of the projected test process depends only on the component of the leading approximation error orthogonal to the parameter derivatives of the approximated conditional probability of $Z>t$. 
We give conditions under which projection eliminates this leading bias term and derive bounds on the remaining Pareto approximation bias. 
Under these conditions, lower thresholds suffice for
asymptotically correct size.

Third, we develop a novel threshold selection criterion specifically for conditional tail index specification tests. 
Existing threshold rules assess the accuracy of tail index estimation \citep{guillouhall2001,schneider2021}, parameter stability across thresholds \citep{wadsworth2016}, or goodness of fit above the threshold \citep{bader2018}. 
Our criterion instead targets the Pareto approximation bias after projection, rather than the entire conditional mean of the test process, which can also reflect misspecification of the conditional tail index model. 
It uses the bias reduction achieved by projection to guide threshold selection, with the objective of limiting size distortions due to the Pareto approximation while retaining more observations for testing. 
In simulations, the projected tests using the new selection criterion generally attain higher power than those using the threshold selection method of \citet{wangtsai2009}.

The remainder of the paper is organized as follows.
Section~\ref{sec:setup} introduces the conditional tail index model.
Section~\ref{sec:method} develops the proposed tests.
Section~\ref{sec:asym_results} establishes asymptotic properties of the tests.
Section~\ref{sec:threshold} presents the threshold selection rule.
Sections~\ref{sec:simulations} and~\ref{sec:real-data} present simulations and a real-data application, respectively.
The Supplementary Material contains proofs, technical details, and additional simulation results.

\section{Conditional tail index models}
\label{sec:setup}

\subsection{Model and null hypothesis}

Let $\{(Y_i,\X_i)\}_{i=1}^{n}$ be independent copies of $(Y,\X)$, where $Y$ is the outcome and $\X=(X_1,\ldots,X_p)\trans$ is a covariate vector taking values in $\mathcal X\subset\R^p$. 
Denote the conditional survival function of $Y$ given $\X=\x$ by $\overline F(y\mid\x)=\Prob(Y>y\mid\X=\x)$. 
For each $\x\in\mathcal X$, we assume that $\overline F(y\mid\x)$ is heavy-tailed in the sense that its upper tail is regularly varying.
Specifically, there exists a function $\gamma_0(\x)>0$ such that
\begin{equation}
 \lim_{y\to\infty}\overline F(ty\mid\x)\big/\overline F(y\mid\x) = t^{-1/\gamma_0(\x)},\quad \text{for all~~} t>0.
 \label{eq:regularvariation}
\end{equation}
Here, $\gamma_0(\x)$ is the conditional tail index at $\x$, with larger values corresponding to heavier upper tails \citep{dehaan2006}.

We are interested in whether a specified parametric model correctly describes how the conditional tail index varies with the covariates.
Let $\{\gamma(\cdot;\tht):\tht\in\Theta\}$ denote the proposed family of positive functions on $\mathcal X$, where $\Theta\subset\R^d$.
The null hypothesis is that $\gamma_0$ belongs to this family, that is,
\begin{equation}
 H_0:\quad \gamma_0(\x)=\gamma(\x;\tht_0)
 \quad\text{for some }\tht_0\in\Theta,
 \label{eq:null}
\end{equation}
where the equality holds for almost every $\x$ with respect to the distribution of $\X$. 
Under the alternative, $\Prob(\gamma_0(\X)\neq\gamma(\X;\tht))>0$ for every $\tht\in\Theta$.
A common specification is the log-linear model
$\gamma(\x;\tht)=\exp[\{\boldsymbol f(\x)\}\trans\tht]$
\citep{wangtsai2009}, where $\boldsymbol f(\x)\in\R^d$
is a specified vector of regressors.

We consider the conditional Hall class of heavy-tailed distributions, which is widely used in the extreme value literature \citep{hall1982,wangtsai2009,li2022,huanglipeng2025}.

\begin{assumption}
\label{ass:hall}
    The conditional tail index satisfies $0<\inf_{\x\in\mathcal X}\gamma_0(\x)\leq\sup_{\x\in\mathcal X}\gamma_0(\x)<\infty$.
    There exist functions $c,\beta,\eta:\mathcal X\to(0,\infty)$, a bounded function $b:\mathcal X\to\R$, and constants $C_0>0$ and $y_0>0$ such that, for every $\x\in\mathcal X$ and every $y\geq y_0$,
    \begin{equation}
     \overline F(y\mid\x) = c(\x)y^{-1/\gamma_0(\x)}
         \{1+b(\x)y^{-\beta(\x)}+R(y,\x)\},
     \label{eq:hall-second-order}
    \end{equation}
    with $|R(y,\x)| \leq C_0 y^{-\beta(\x)-\eta(\x)}$.
    In addition, $\underline\beta:=\inf_{\x\in\mathcal X}\beta(\x)>0$ and $\underline\eta:=\inf_{\x\in\mathcal X}\eta(\x)>0$.
\end{assumption}

The Hall class includes many familiar heavy-tailed distributions, such as the Burr, Fr\'echet, and Student-$t$. 
An exact Pareto tail at $\x$ is included with $b(\x)=0$ and $R(y,\x)=0$ for all $y\geq y_0$ in \eqref{eq:hall-second-order}. 

For the theoretical analysis, let $u>0$ be a deterministic threshold.
Denote by $N=\sum_{i=1}^{n}\ind\{Y_i>u\}$ the number of observations with $Y_i>u$, and index these observations by $i=1,\ldots,N$. 
For the analysis under $H_0$, we condition on
$\mathcal C_N=\sigma(N,\X_1,\ldots,\X_N)$.
We consider sequences along which $u\to\infty$ and $N\to\infty$ as $n\to\infty$.

The conditional survival function of $Z=\log(Y/u)$ for $Y>u$ is
\begin{equation*}
    p_t(\x) = \Prob(Z>t\mid Y>u,\X=\x) = \overline F(ue^t\mid\x)\big/\overline F(u\mid\x),\quad \text{for all~~} t\geq0.
\end{equation*}
By \eqref{eq:regularvariation}, $p_t(\x)$ converges to
$\exp\{-t/\gamma_0(\x)\}$ as $u\to\infty$ for each
$\x\in\mathcal X$ and each fixed $t\geq0$.
Thus, the conditional distribution of the log excess is approximated by an exponential distribution with mean $\gamma_0(\x)$. 
This approximation provides the basis for estimating the null parameter and constructing the specification tests.
Let $Z_i=\log(Y_i/u)$ for $i=1,\ldots,N$. 
Conditional on $\mathcal C_N$, the log excesses are independent and satisfy $\Prob(Z_i>t\mid\mathcal C_N)=p_t(\X_i)$.

\subsection{Estimation of the null parameter}
\label{sec:estimation}

We estimate the null parameter $\tht_0$ by fitting the exponential approximation to the log excesses. 
Let $\pn$ denote the empirical average over $(Y_1,\X_1),\ldots,(Y_N,\X_N)$. 
For a function $f$, we write
$\pn\{f(Y,\X)\}=N^{-1}\sum_{i=1}^N f(Y_i,\X_i)$. 
The average negative log quasi-likelihood is
\begin{equation*}
    Q_N(\tht) = \pn\left\{\log\gamma(\X;\tht)+Z/\gamma(\X;\tht)\right\}.
\end{equation*}
Let $\widehat\tht$ and $\tht_u$ denote the local minimizers of $Q_N(\tht)$ and $\E\{Q_N(\tht)\mid\mathcal C_N\}$, respectively, in a neighborhood of $\tht_0$. 
Here $\tht_u$ is the pseudo-true parameter at threshold $u$ \citep{white1982}. %

Write $\mu_{u,i}=\E(Z_i\mid\mathcal C_N)$ as the conditional mean of the log excess at threshold $u$ and $\vv(\x;\tht)=\partial_{\tht}\log\gamma(\x;\tht)$.
Then $\tht_u$ satisfies the score equation
\begin{equation}
 N^{-1}\sum_{i=1}^N
 \vv(\X_i;\tht_u)
 \left\{\mu_{u,i}\big/\gamma(\X_i;\tht_u)-1\right\}
 =\boldsymbol 0.
 \label{eq:pseudotrue}
\end{equation}
Under $H_0$, the model specifies the conditional tail index $\gamma_0(\X_i)$, whereas the score equation depends on $\mu_{u,i}$.
If the conditional tails are exactly Pareto above $u$, then
$\mu_{u,i}=\gamma_0(\X_i)$ for every $i$, and local identification gives $\tht_u=\tht_0$. 
More generally, Assumption~\ref{ass:hall} implies $\max_{1\leq i\leq N}|\mu_{u,i}/\gamma_0(\X_i)-1|=O(u^{-\underline\beta})$, so the pseudo-true parameter $\tht_u$ can differ from $\tht_0$ even under the null hypothesis.

Define
$\Sigma_{\X,0}=N^{-1}\sum_{i=1}^N\vv(\X_i;\tht_0)\{\vv(\X_i;\tht_0)\}\trans$ and let $\lambda_{\min}(\Sigma_{\X,0})$ be the smallest eigenvalue of $\Sigma_{\X,0}$.

\begin{assumption}
\label{ass:regular}
    (i) The parameter $\tht_0$ is an interior point of $\Theta$. In a neighborhood of $\tht_0$, $\gamma(\x;\tht)$ is bounded above and away from zero uniformly over $\x\in\mathcal X$ and is twice continuously differentiable with respect to $\tht$. Its first two parameter derivatives are uniformly bounded over $\x\in\mathcal X$ and this neighborhood, and its second derivative is Lipschitz continuous in $\tht$ uniformly over $\x\in\mathcal X$. 
    (ii) There exists a constant $c_0>0$ such that $\lambda_{\min}(\Sigma_{\X,0})\geq c_0$ for all sufficiently large $N$.
\end{assumption}

For compact $\mathcal X$, joint continuity of $\gamma(\x;\tht)$ and its first two derivatives in $(\x,\tht)$ ensures the local boundedness and continuity requirements in Assumption~\ref{ass:regular}(i). 
A uniformly bounded third derivative with respect to $\tht$ is sufficient for the Lipschitz condition. 
Part~(ii) of Assumption~\ref{ass:regular} requires the observations with $Y_i>u$ to remain informative about every parameter direction. 
Under a common threshold, their covariates may concentrate near values with the heaviest tails, so this condition can fail \citep{yangrank2026}.
Section~\ref{sec:scope} studies the validity of the proposed tests when $\lambda_{\min}(\Sigma_{\X,0})\to0$.

\begin{theorem}
\label{thm:hall-estimation}
    Suppose that $H_0$ and Assumptions~\ref{ass:hall} and~\ref{ass:regular} hold.
    For all sufficiently large $u$ and $N$, the function $\tht\mapsto\E\{Q_N(\tht)\mid\mathcal C_N\}$ has a unique local minimizer $\tht_u$ in a sufficiently small neighborhood of $\tht_0$.
    With probability tending to one, $Q_N(\tht)$ has a unique local minimizer $\widehat\tht$ in the same neighborhood.
    In addition, we have 
    
    \noindent (i) The discrepancy between $\tht_u$ and $\tht_0$ satisfies
    $\tht_u-\tht_0=\boldsymbol\delta_{u,N}+O(u^{-2\underline\beta}+u^{-\underline\beta-\underline\eta})$, where
    \begin{equation*}
        \boldsymbol\delta_{u,N} = -N^{-1}\Sigma_{\X,0}^{-1}\sum_{i=1}^N\vv(\X_i;\tht_0)\,b(\X_i)u^{-\beta(\X_i)}\frac{\beta(\X_i)\gamma_0(\X_i)}{1+\beta(\X_i)\gamma_0(\X_i)}.
    \end{equation*}

    \noindent (ii) $N^{1/2}\Sigma_{\X,0}^{1/2}(\widehat\tht-\tht_u) = N^{-1/2}\Sigma_{\X,0}^{-1/2}\sum_{i=1}^N\vv(\X_i;\tht_0)(Z_i-\mu_{u,i})/\gamma_0(\X_i)+\Op(N^{-1/2}+u^{-\underline\beta})$ and
    \begin{equation}
        N^{1/2}\Sigma_{\X,0}^{1/2}(\widehat\tht-\tht_u) \xrightarrow{d} N(\boldsymbol 0,\boldsymbol{\mathsf I}_d).
        \label{eq:hall-centered-clt}
    \end{equation}
    If, in addition, $N^{1/2}(u^{-2\underline\beta}+u^{-\underline\beta-\underline\eta})\to0$, then
    \begin{equation}
        N^{1/2}\Sigma_{\X,0}^{1/2}\bigl(\widehat\tht-\tht_0-\boldsymbol\delta_{u,N}\bigr) \xrightarrow{d} N(\boldsymbol 0,\boldsymbol{\mathsf I}_d).
        \label{eq:hallbiasednormal}
    \end{equation}
\end{theorem}

Theorem~\ref{thm:hall-estimation} separates the estimation error $\widehat\tht-\tht_u=\Op(N^{-1/2})$ from the difference $\tht_u-\tht_0=O(u^{-\underline\beta})$ induced by the Pareto approximation. 
The limit in \eqref{eq:hall-centered-clt} remains valid when this difference is nonnegligible relative to $N^{-1/2}$. 
The additional rate condition in part~(ii) ensures that, after subtracting the leading bias $\boldsymbol\delta_{u,N}$, the remaining bias is negligible. 

\section{Neyman-orthogonalized specification tests}
\label{sec:method}

\subsection{Test processes}

The exponential approximation in Section~\ref{sec:setup} relates the conditional tail index to the probability that an observation with $Y>u$ also exceeds the higher level $ue^t$.
For $\tht\in\Theta$ and $t\geq0$, let
$S_t(\x;\tht)=\exp\{-t/\gamma(\x;\tht)\}$
denote this conditional probability under the exponential approximation.
For each observation with $Y_i>u$, the indicator
$\ind\{Z_i>t\}=\ind\{Y_i>ue^t\}$ records whether the higher level is exceeded.
Since
$\E(\ind\{Z_i>t\}\mid\mathcal C_N)=p_t(\X_i)$,
a natural way to assess the null specification is to compare these
indicators with the probabilities implied by the exponential approximation.
We therefore define the residual
\begin{equation}
    r_i(t;\tht) = \ind\{Z_i>t\}-S_t(\X_i;\tht),
    \qquad t\geq0.
    \label{eq:residual}
\end{equation}
Under $H_0$, the conditional mean of the residual evaluated at the true parameter is
$\E\{r_i(t;\tht_0)\mid\mathcal C_N\}=p_t(\X_i)-S_t(\X_i;\tht_0)$. 
This mean is zero if the conditional tail is exactly Pareto above $u$. 
More generally, even when the conditional tail index model is correctly specified, the residual can have a nonzero conditional mean because the Pareto approximation need not be exact.

A misspecified model may overstate the conditional probability at some covariate values and understate it at others.
Such discrepancies can cancel when the residuals are averaged over all observations with $Y_i>u$.
Following \citet{bierens1990} and \citet{stute1997}, we instead accumulate the residuals over covariate regions of the form $\{\X_i\leq\x\}$.
Let $h_{\x}(\boldsymbol z)=\ind\{\boldsymbol z\leq\x\}$, where the
inequality is coordinatewise, and let $\mathcal T$ be a compact interval in $(0,\infty)$ or a finite set of positive log excess levels.
Using the estimator $\widehat\tht$ from  Section~\ref{sec:estimation}, for $t\in\mathcal T$ and $\x\in\R^p$, define the unprojected test process
\begin{equation}
    U_N^0(t,\x) = N^{-1/2}\sum_{i=1}^N h_{\x}(\X_i)\,r_i(t;\widehat\tht).
    \label{eq:rawprocess}
\end{equation}
Varying $\x$ allows discrepancies across covariate regions to be examined without nonparametric estimation of the conditional survival function.

To examine the effect of parameter estimation, write
\begin{equation}
    \dd_t(\x;\tht) 
    = \partial_{\tht}S_t(\x;\tht) 
    = \frac{t}{\gamma(\x;\tht)}S_t(\x;\tht)\vv(\x;\tht).
    \label{eq:derivative}
\end{equation}
Under $H_0$ and Assumptions~\ref{ass:hall} and~\ref{ass:regular}, a Taylor expansion about $\tht_u$, together with Theorem~\ref{thm:hall-estimation}, gives
\begin{equation}
    \begin{aligned}
        U_N^0(t,\x) ={}& N^{-1/2}\sum_{i=1}^N h_{\x}(\X_i)\{\ind\{Z_i>t\}-p_t(\X_i)\} -\left[\pn\{h_{\x}(\X)\dd_t(\X;\tht_u)\}\right]\trans N^{1/2}(\widehat\tht-\tht_u)\\
        &+N^{1/2}\pn\left[h_{\x}(\X)\{p_t(\X)-S_t(\X;\tht_u)\}\right]+\Op(N^{-1/2}),
    \end{aligned}
    \label{eq:raw-fitted-expansion}    
\end{equation}
uniformly over $t\in\mathcal T$ and $\x\in\R^p$. 
The first summation term has zero mean conditional on $\mathcal C_N$ and is uniformly $\Op(1)$. 
As $\widehat\tht-\tht_u=\Op(N^{-1/2})$, the second term is also uniformly $\Op(1)$ and generally contributes to the asymptotic distribution. 
The third term is the conditional mean of the unprojected test process with the parameter fixed at $\tht_u$ and is uniformly $O(N^{1/2}u^{-\underline\beta})$. 
Thus parameter estimation can affect the asymptotic distribution of $U_N^0(t,\x)$ even when the exponential approximation is exact and the third term vanishes. 
A multiplier bootstrap that reweights the summands in \eqref{eq:rawprocess} while holding $\widehat\tht$ fixed generally does not reproduce this parameter estimation contribution.

The second term in \eqref{eq:raw-fitted-expansion} suggests replacing $h_{\x}(\X_i)$ by weights orthogonal to the components of $\dd_t(\X_i;\tht)$. 
We use the least-squares residuals from regressing $h_{\x}(\X_i)$ on these components across $i=1,\ldots,N$. 
When the matrix with rows $\{\dd_t(\X_i;\tht)\}\trans$ has full column rank, the resulting weight is
\begin{equation}
 h_{t,\x}^{\perp}(\X_i;\tht)
 =h_{\x}(\X_i)
  -\{\dd_t(\X_i;\tht)\}\trans
   \bigl(\pn[\dd_t(\X;\tht)\{\dd_t(\X;\tht)\}\trans]\bigr)^{-1}
   \pn\{\dd_t(\X;\tht) h_{\x}(\X)\}.
 \label{eq:weightprojection}
\end{equation}
The least-squares normal equations give the sample orthogonality condition
\begin{equation}
    \pn\left\{h_{t,\x}^{\perp}(\X;\tht)\dd_t(\X;\tht)\right\}
    =\boldsymbol 0.
    \label{eq:exactorthogonality}
\end{equation}
Replacing the region indicators in \eqref{eq:rawprocess} by these weights evaluated at $\widehat\tht$ gives the Neyman-orthogonalized test process \citep{neyman1959,chernozhukov2018},
\begin{equation*}
    U_N(t,\x)
    = N^{-1/2}\sum_{i=1}^N h_{t,\x}^{\perp}(\X_i;\widehat\tht)r_i(t;\widehat\tht).
\end{equation*}

To verify Neyman orthogonality, let $\E_0$ denote expectation with $\X_1,\ldots,\X_N$ held fixed and the log excesses independently distributed as exponential random variables with means $\gamma(\X_i;\tht_0)$. 
Keeping this distribution fixed, we have, for every $\boldsymbol a\in\R^d$, 
\begin{equation}
 \left.
 \frac{d}{d\varepsilon}\,
 \E_0\!\left[
   N^{-1}\sum_{i=1}^N
   h_{t,\x}^{\perp}(\X_i;\tht_0+\varepsilon\boldsymbol a)\,
   r_i(t;\tht_0+\varepsilon\boldsymbol a)
 \right]
 \right|_{\varepsilon=0}
 =0.
 \notag
\end{equation}
The term from differentiating the weights vanishes because $\E_0\{r_i(t;\tht_0)\}=0$ for every $i$, and the term from differentiating the residuals vanishes by \eqref{eq:exactorthogonality}.
Under the actual conditional distribution of the log excess, \eqref{eq:exactorthogonality} still holds, but $p_t(\X_i)-S_t(\X_i;\tht_0)$ need not be zero even under $H_0$.
The contribution from differentiating the weights can therefore remain.
Section~\ref{sec:mechanism} quantifies this contribution and characterizes the Pareto approximation bias after projection.

\subsection{Test statistics and multiplier bootstrap}
\label{sec:statistic}\label{sec:bootstrap}

Let $F_{\X,N}$ denote the empirical distribution of $\X_1,\ldots,\X_N$ and let $W$ be a probability measure on $\mathcal T$.
We construct the CvM and KS statistics as
\begin{equation*}
    T_N^{\mathrm{CvM}}=\int_{\mathcal T}\int_{\R^p} \{U_N(t,\x)\}^2\,dF_{\X,N}(\x)\,dW(t) 
    \qquad \text{~~and~~} \qquad
    T_N^{\mathrm{KS}}=\sup_{t\in\mathcal T,\,\x\in\R^p}|U_N(t,\x)|.
\end{equation*}
The CvM statistic integrates squared discrepancies over the covariates and log excess levels, whereas the KS statistic takes the largest absolute discrepancy over all $(t,\x)\in\mathcal T\times\R^p$ \citep{stute1997,santanna2019}. 

We apply a multiplier bootstrap \citep{kojadinovicyan2012,escancianogoh2014,santanna2019} to approximate the distributions of $T_N^{\mathrm{CvM}}$ and $T_N^{\mathrm{KS}}$, and to obtain the corresponding critical values. 
Let $\xi_1,\ldots,\xi_N$ be i.i.d.\ multipliers, independent of the data, with $\E(\xi_1)=0$, $\E(\xi_1^2)=1$, and $\E|\xi_1|^{2+\delta}<\infty$ for some $\delta>0$.
Define the multiplier bootstrap process by
\begin{equation*}
    U_N^*(t,\x) = N^{-1/2}\sum_{i=1}^N\xi_i\,h_{t,\x}^{\perp}(\X_i;\widehat\tht)\,r_i(t;\widehat\tht).
\end{equation*}

The bootstrap counterparts of the test statistics are
\begin{equation*}
    T_N^{*,\mathrm{CvM}} = \int_{\mathcal T}\int_{\R^p} \{U_N^*(t,\x)\}^2\,dF_{\X,N}(\x)\,dW(t) 
    \qquad \text{~~and~~} \qquad
    T_N^{*,\mathrm{KS}} = \sup_{t\in\mathcal T,\,\x\in\R^p}|U_N^*(t,\x)|.
\end{equation*}
For either test, write $T_N$ for the projected test statistic and $T_N^{*(b)}$ for its counterpart in bootstrap draw $b$.
The bootstrap $p$-value based on $B$ independent draws is
\begin{equation}
 \widehat p
 =\left[1+\sum_{b=1}^B\ind\{T_N^{*(b)}\geq T_N\}\right]/(B+1).
 \label{eq:pvalue}
\end{equation}
At a prespecified nominal level $\alpha\in(0,1)$, we reject the null hypothesis if the bootstrap $p$-value for the chosen statistic satisfies $\widehat p\leq\alpha$.

\subsection{Effects of projection on the test process}
\label{sec:mechanism}

We now expand $U_N(t,\x)$ about $\tht_u$ and compare it with the expansion of $U_N^0(t,\x)$ in \eqref{eq:raw-fitted-expansion}. 
Write $e_i(t)=p_t(\X_i)-S_t(\X_i;\tht_u)$ for the conditional mean of $r_i(t;\tht_u)$, and define
\begin{equation*}
    \mu_N(t,\x) = N^{1/2}\pn\left[h_{t,\x}^{\perp}(\X;\tht_u)\{p_t(\X)-S_t(\X;\tht_u)\}\right] = N^{-1/2}\sum_{i=1}^{N}h_{t,\x}^{\perp}(\X_i;\tht_u)e_i(t).
\end{equation*}
Thus $\mu_N(t,\x)$ is the conditional mean of the projected test process with the parameter fixed at $\tht_u$.
Under $H_0$, it represents the Pareto approximation bias after projection. %

\begin{proposition}
\label{prop:fitted-expansions}
    Suppose that $H_0$ and Assumptions~\ref{ass:hall} and~\ref{ass:regular} hold. 
    Then the following expansion holds uniformly over $t\in\mathcal T$ and $\x\in\R^p$:
    \begin{equation}
        U_N(t,\x) ={} N^{-1/2}\sum_{i=1}^N h_{t,\x}^{\perp}(\X_i;\tht_u)\{\ind\{Z_i>t\}-p_t(\X_i)\} +\mu_N(t,\x)+\Op\left(N^{-1/2}+u^{-\underline\beta}\right).
        \label{eq:wp-fitted-expansion}
    \end{equation}
    The first summation term in \eqref{eq:wp-fitted-expansion} has zero mean conditional on $\mathcal C_N$ and is uniformly $\Op(1)$.
\end{proposition}

The term arising from differentiation of $S_t(\X;\tht)$ in \eqref{eq:raw-fitted-expansion} vanishes by the orthogonality condition \eqref{eq:exactorthogonality}. 
Replacing $\widehat\tht$ with $\tht_u$ in the projected weight contributes $\Op(N^{-1/2})$ to the first summation term in \eqref{eq:wp-fitted-expansion} and $\Op(u^{-\underline\beta})$ to the term involving the conditional residual means $e_i(t)$. 
Together with the remainder of order $\Op(N^{-1/2})$, these contributions give the remainder in \eqref{eq:wp-fitted-expansion}.
Consequently, $U_N(t,\x)$ is asymptotically equivalent to the test process with the parameter fixed at $\tht_u$.

The conditional mean term $\mu_N(t,\x)$ remains in the expansion of the projected test process and can be nonnegligible even under $H_0$.
It depends on both the error in the Pareto approximation at $\tht_0$ and the parameter shift $\tht_u-\tht_0$.
We examine how projection acts on these two contributions.
Under Assumption~\ref{ass:hall}, the leading term in $p_t(\x)-S_t(\x;\tht_0)$ is
\begin{equation}
    a_{u,t}(\x)
    =S_t(\x;\tht_0)b(\x)u^{-\beta(\x)}
    \{e^{-\beta(\x)t}-1\}.
    \label{eq:leading-survivor-correction}
\end{equation}
Theorem~\ref{thm:hall-estimation} gives $\tht_u-\tht_0=O(u^{-\underline\beta})$. 
A Taylor expansion of $S_t(\X_i;\tht_u)$ about $\tht_0$ therefore yields
\begin{equation*}
    \notag
    e_i(t) = a_{u,t}(\X_i)-\{\dd_t(\X_i;\tht_0)\}\trans(\tht_u-\tht_0)
    +O\left(u^{-2\underline\beta}+u^{-\underline\beta-\underline\eta}\right),
\end{equation*}
uniformly over $i=1,\ldots,N$ and $t\in\mathcal T$. 
The term linear in $\tht_u-\tht_0$ is removed by the projection normal equations at $\tht_0$. 
Moreover, replacing $\tht_u$ by $\tht_0$ in the projection weights changes $\mu_N(t,\x)$ by $O(N^{1/2}u^{-2\underline\beta})$.
These observations yield the following expansions.

\begin{proposition}
\label{prop:projected-bias-expansion}
    Suppose that $H_0$ and Assumptions~\ref{ass:hall} and~\ref{ass:regular} hold. 
    Then the following expansions hold uniformly over $t\in\mathcal T$ and $\x\in\R^p$:
    \begin{equation}
        \mu_N(t,\x) = N^{1/2}\pn\left\{h_{t,\x}^{\perp}(\X;\tht_0)a_{u,t}(\X)\right\}+O\left(N^{1/2}(u^{-2\underline\beta}+u^{-\underline\beta-\underline\eta})\right).
        \label{eq:wp-leading-mean}
    \end{equation}
    and 
    \begin{equation}
        \begin{aligned}
            U_N(t,\x)
            ={}&N^{-1/2}\sum_{i=1}^N
            h_{t,\x}^{\perp}(\X_i;\tht_u)
            \{\ind\{Z_i>t\}-p_t(\X_i)\}\\
            &+N^{1/2}\pn\left\{
            h_{t,\x}^{\perp}(\X;\tht_0)a_{u,t}(\X)
            \right\}+\Op\left(
            N^{-1/2}+N^{1/2}(u^{-2\underline\beta}+u^{-\underline\beta-\underline\eta})\right).
        \end{aligned}
        \label{eq:combined-process-rate}
    \end{equation}
\end{proposition}

Under $H_0$, the term $N^{1/2}\pn\{h_{t,\x}^{\perp}(\X;\tht_0)a_{u,t}(\X)\}$ in \eqref{eq:combined-process-rate} is the leading term of the Pareto approximation bias in the projected test process.
The full bias is $\mu_N(t,\x)$, whose expansion in \eqref{eq:wp-leading-mean} also includes a remainder of order $N^{1/2}(u^{-2\underline\beta}+u^{-\underline\beta-\underline\eta})$. 
Since projection does not increase the empirical second moment of the indicator weights, the Cauchy--Schwarz inequality gives
\begin{equation*}
    \sup_{t\in\mathcal T,\,\x\in\R^p}
    \left|
    \pn\left\{
    h_{t,\x}^{\perp}(\X;\tht_0)a_{u,t}(\X)
    \right\}
    \right|
    \leq
    \sup_{t\in\mathcal T}
    \left(\pn\left[\{a_{u,t}(\X)\}^2\right]\right)^{1/2}
    =O(u^{-\underline\beta}).
\end{equation*}
Thus $N^{1/2}u^{-\underline\beta}\to0$ is sufficient for this
contribution to be negligible. 

Substitution from \eqref{eq:weightprojection} gives
\begin{equation}
    \begin{aligned}
        \pn\left\{h_{t,\x}^{\perp}(\X;\tht_0)a_{u,t}(\X)\right\}
         = & \pn\{h_{\x}(\X)a_{u,t}(\X)\} \\
        & - \left[\pn\{h_{\x}(\X)\dd_t(\X;\tht_0)\}\right]\trans
        \left(\pn[\dd_t(\X;\tht_0)\{\dd_t(\X;\tht_0)\}\trans]\right)^{-1}\pn\{\dd_t(\X;\tht_0)a_{u,t}(\X)\}.
    \end{aligned}
    \label{eq:projectedbiasidentity}
\end{equation}
Equivalently, it is the average over the region indexed by $\x$ of the residuals from regressing $a_{u,t}(\X_i)$ on $\dd_t(\X_i;\tht_0)$.
Only the component of $a_{u,t}$ orthogonal to the components of $\dd_t(\cdot;\tht_0)$ contributes to the leading term of $\mu_N(t,\x)$ in \eqref{eq:wp-leading-mean}.

A sufficient condition for $\pn\{h_{t,\x}^{\perp}(\X;\tht_0)a_{u,t}(\X)\}$ to vanish is that, for every $t\in\mathcal T$ and all sufficiently large $u$, there exists $\boldsymbol c_{u,t}\in\R^d$ such that
\begin{equation}
    \frac{\gamma_0(\x)b(\x)}{t}
    u^{-\beta(\x)}
    \{e^{-\beta(\x)t}-1\}
    =
    \{\vv(\x;\tht_0)\}\trans\boldsymbol c_{u,t},
    \qquad \x\in\mathcal X.
    \label{eq:hall-alignment}
\end{equation}
Since $\dd_t(\x;\tht_0)=\{t/\gamma_0(\x)\}S_t(\x;\tht_0)\vv(\x;\tht_0)$, condition \eqref{eq:hall-alignment} is equivalent to
\begin{equation}
    \notag
    a_{u,t}(\x)
    =\{\dd_t(\x;\tht_0)\}\trans\boldsymbol c_{u,t},
    \qquad \x\in\mathcal X.
\end{equation}
The projection normal equations therefore imply the following result.

\begin{corollary}
\label{cor:projection-rate}
    Under the conditions of Proposition~\ref{prop:projected-bias-expansion}, suppose that \eqref{eq:hall-alignment} holds.
    Then the leading term in $\mu_N(t,\x)$ is identically zero, and
    \begin{equation}
        \sup_{t\in\mathcal T,\,\x\in\R^p}
        |\mu_N(t,\x)|
        =
        O\left(
        N^{1/2}(
        u^{-2\underline\beta}
        +u^{-\underline\beta-\underline\eta})
        \right).
        \label{eq:exact-projected-mean-rate}
    \end{equation}
\end{corollary}

Corollary~\ref{cor:projection-rate} shows that, under \eqref{eq:hall-alignment}, projection can reduce the order of the Pareto approximation bias in the test process.
The general sufficient condition $N^{1/2}u^{-\underline\beta}\to0$ is then replaced by
\begin{equation*}
    N^{1/2}(
    u^{-2\underline\beta}
    +u^{-\underline\beta-\underline\eta}
    )\to0.
\end{equation*}
In particular, if $\underline\eta\geq\underline\beta$, the sufficient condition becomes $N^{1/2}u^{-2\underline\beta}\to0$.
This requires $u\gg N^{1/(4\underline\beta)}$, compared with $u\gg N^{1/(2\underline\beta)}$ under the general bound.
Thus a lower threshold can suffice for the mean of the projected test process to be negligible, while $\tht_u-\tht_0$ can remain of order $u^{-\underline\beta}$.

\begin{remark}
    When $\beta(\x)\equiv\beta$ for a constant $\beta$, the expansion in Proposition~\ref{prop:projected-bias-expansion} becomes
    \begin{equation}
        \mu_N(t,\x) 
        = N^{1/2}u^{-\beta}(e^{-\beta t}-1)\pn\left\{h_{t,\x}^{\perp}(\X;\tht_0)S_t(\X;\tht_0)b(\X)\right\}
        + O\left(N^{1/2}(u^{-2\beta}+u^{-\beta-\underline\eta})\right),
        \label{eq:common-beta-bias}
    \end{equation}
    uniformly over $t\in\mathcal T$ and $\x\in\R^p$. 
    For any fixed $t\in\mathcal T$, the remainder in \eqref{eq:common-beta-bias} is $o(N^{1/2}u^{-\beta})$. 
    Thus, if the empirical average in \eqref{eq:common-beta-bias} converges to a nonzero constant at some fixed $(t,\x)$, the Pareto approximation bias at that point remains of order $N^{1/2}u^{-\beta}$.
    For a common $\beta$, condition~\eqref{eq:hall-alignment}
    reduces to the existence of a vector $\boldsymbol c\in\R^d$
    such that
    $\gamma_0(\x)b(\x) = \{\vv(\x;\tht_0)\}\trans\boldsymbol c$
    for all $\x\in\mathcal X$.
    In a log-linear model, this condition means
    that $\gamma_0(\x) b(\x)$ is a linear combination of the regressors.
\end{remark}

A reduction in the error order does not require the exact representation \eqref{eq:hall-alignment}.
More generally, suppose that, for some $\kappa\geq\underline\beta$,
\begin{equation}
    \sup_{t\in\mathcal T,\,\x\in\R^p}
    \left|
    \pn\left\{
    h_{t,\x}^{\perp}(\X;\tht_0)a_{u,t}(\X)
    \right\}
    \right|
    =O(u^{-\kappa}).
    \label{eq:regional-bias-rate}
\end{equation}
For example, this condition holds if
$a_{u,t}(\X_i)=\{\dd_t(\X_i;\tht_0)\}\trans\boldsymbol c_{u,t}+\varepsilon_{u,t,i}$ and
$\sup_{t\in\mathcal T}\bigl\{N^{-1}\sum_{i=1}^N\varepsilon_{u,t,i}^2\bigr\}^{1/2}=O(u^{-\kappa})$.
Under \eqref{eq:regional-bias-rate}, let
$q=\min\{\kappa,\,2\underline\beta,\,\underline\beta+\underline\eta\}$. 
Proposition~\ref{prop:projected-bias-expansion} then gives
$\sup_{t\in\mathcal T,\,\x\in\R^p}|\mu_N(t,\x)|=O(N^{1/2}u^{-q})$,
so $N^{1/2}u^{-q}\to0$ is sufficient for both the mean and the remainder in \eqref{eq:combined-process-rate} to vanish.
Under \eqref{eq:hall-alignment}, the same conclusion holds with $q=\min\{2\underline\beta,\underline\beta+\underline\eta\}$.

The following example illustrates a model in which projection eliminates the leading Pareto approximation bias.
Further examples are given in Section~A.2 of the Supplementary Material.

\begin{example}
    Consider the conditional Burr model with 
    $\overline F(y\mid\x)=(1+y^2)^{-1/\{2\gamma_0(\x)\}}$ for $y>0$, where $\log\gamma_0(\x)=\{\boldsymbol f(\x)\}\trans\tht_0$ and the first component of $\boldsymbol f(\x)$ is one.
    Here $b(\x)=-1/\{2\gamma_0(\x)\}$ and $\beta(\x)=\eta(\x)=2$.
    Since $\vv(\x;\tht_0)=\boldsymbol f(\x)$, condition \eqref{eq:hall-alignment} holds with
    $\boldsymbol c_{u,t}
        =
        (1-e^{-2t})({2tu^2})^{-1}
        (1,0,\ldots,0)\trans$.
    Thus projection removes the leading Pareto approximation bias even when the tail index varies with the covariates, and $\mu_N(t,\x)$ is uniformly $O(N^{1/2}u^{-4})$.
    The sufficient condition is therefore $N^{1/2}u^{-4}\to0$, compared with $N^{1/2}u^{-2}\to0$ obtained from the general approximation bound.
\end{example}

Under $H_0$, Proposition~\ref{prop:projected-bias-expansion} and Corollary~\ref{cor:projection-rate} relate the sufficient threshold conditions to the Pareto approximation bias after projection. 
Section~\ref{sec:threshold} proposes a novel threshold selection rule that estimates this bias and selects the lowest candidate satisfying a tolerance criterion.

\section{Asymptotic theory}\label{sec:asym_results}

\subsection{Asymptotic size}
\label{sec:theory}

We first establish Gaussian and multiplier approximations for the CvM and KS statistics based on the centered test process 
$U_N(t,\x)-\mu_N(t,\x)$ under $H_0$. 
We then give conditions under which the Pareto approximation bias $\mu_N(t,\x)$ is negligible and the bootstrap tests have asymptotic size $\alpha$. 
Proposition~\ref{prop:fitted-expansions} gives
\begin{equation}
    \sup_{t\in\mathcal T,\,\x\in\R^p}
        \left|U_N(t,\x)-\mu_N(t,\x)-N^{-1/2}\sum_{i=1}^N h_{t,\x}^{\perp}(\X_i;\tht_u)\{\ind\{Z_i>t\}-p_t(\X_i)\}\right|
        = \Op\left(N^{-1/2}+u^{-\underline\beta}\right).
    \label{eq:uniform-wp-expansion}
\end{equation}
The weighted empirical process $N^{-1/2}\sum_{i=1}^N h_{t,\x}^{\perp}(\X_i;\tht_u)\{\ind\{Z_i>t\}-p_t(\X_i)\}$ in
\eqref{eq:uniform-wp-expansion} has mean zero conditional on
$\mathcal C_N$ and conditional covariance function
\begin{equation*}
    \Gamma_N(t,\x;s,\boldsymbol y) = N^{-1}\sum_{i=1}^N h_{t,\x}^{\perp}(\X_i;\tht_u)h_{s,\boldsymbol y}^{\perp}(\X_i;\tht_u)\{p_{t\vee s}(\X_i)-p_t(\X_i)p_s(\X_i)\}.
\end{equation*}
Let $G_N=\{G_N(t,\x):(t,\x)\in\mathcal T\times\R^p\}$ be a  Gaussian process with mean zero and covariance function $\Gamma_N$ conditional on $\mathcal C_N$.
Let $\nu_u$ denote the conditional distribution of $\X$ given $Y>u$. 
As this distribution may change with $u$, we formulate the Gaussian approximation using the sequence $\{G_N\}$ without requiring $\Gamma_N$ to converge to a fixed covariance function.

Define the Gaussian counterparts of the test statistics by
$T_{N,G}^{\mathrm{CvM}}=\int_{\mathcal T}\int_{\R^p}\{G_N(t,\x)\}^2\,dF_{\X,N}(\x)\,dW(t)$
and
$T_{N,G}^{\mathrm{KS}}=\sup_{t\in\mathcal T,\,\x\in\R^p}|G_N(t,\x)|$.
The corresponding statistics based on $U_N(t,\x)-\mu_N(t,\x)$ are
$T_{N,c}^{\mathrm{CvM}}=\int_{\mathcal T}\int_{\R^p}\{U_N(t,\x)-\mu_N(t,\x)\}^2\,dF_{\X,N}(\x)\,dW(t)$ and
$T_{N,c}^{\mathrm{KS}}=\sup_{t\in\mathcal T,\,\x\in\R^p}|U_N(t,\x)-\mu_N(t,\x)|$.
For either CvM or KS test, write $T_{N,G}$ and $T_{N,c}$ for the corresponding statistics based on $G_N(t,\x)$ and $U_N(t,\x)-\mu_N(t,\x)$ defined above, respectively. 
We impose the following condition on $\Gamma_N$. 

\begin{assumption}
\label{ass:variation}
    There exists a constant $c_1>0$ such that, for all sufficiently large $N$, 
    $\int_{\mathcal T}\int_{\R^p}\Gamma_N(t,\x;t,\x)\,dF_{\X,N}(\x)\,dW(t)\geq c_1$ for CvM, and
    $\sup_{t\in\mathcal T,\,\x\in\R^p}\Gamma_N(t,\x;t,\x)\geq c_1$ for KS.
\end{assumption}

Assumption~\ref{ass:variation} requires the integrated conditional variance for CvM, or the supremum of the conditional variance for KS, to be bounded away from zero.
It concerns the variation remaining in the test process after projection, whereas Assumption~\ref{ass:regular}(ii) concerns information for estimating $\tht$.
For a scalar covariate, Lemma~A.3 in Section~A.3 of the Supplementary Material verifies both requirements of Assumption~\ref{ass:variation} under $H_0$, without imposing Assumption~\ref{ass:regular}(ii). 
Let $\Prob^*$ and $\E^*$ denote probability and expectation over the multipliers conditional on the observed data.

\begin{proposition}
\label{prop:bounded-drift}
    Under $H_0$ and Assumptions~\ref{ass:hall}, \ref{ass:regular}, and~\ref{ass:variation},
    \begin{equation}
        \begin{aligned}
            & \sup_{z\in\R}\left|\Prob(T_{N,c}\leq z\mid\mathcal C_N)-\Prob(T_{N,G}\leq z\mid\mathcal C_N)\right|\to0 \text{~~and~~} \\
            & \sup_{z\in\R}\left|\Prob^*(T_N^*\leq z)-\Prob(T_{N,G}\leq z\mid\mathcal C_N)\right|\overset{p}{\to}0.
        \end{aligned}
        \label{eq:wp-triangular-bootstrap}
    \end{equation}
\end{proposition}

Proposition~\ref{prop:bounded-drift} shows that the multiplier bootstrap consistently approximates the distributions of the CvM and KS statistics based on $U_N(t,\x)-\mu_N(t,\x)$, even when $\mu_N(t,\x)$ is nonnegligible.
To establish asymptotic size for the bootstrap tests, we impose the following condition on $\mu_N(t,\x)$.

\begin{assumption}
\label{ass:bias}
    For CvM,
    $\int_{\mathcal T}\int_{\R^p}\{\mu_N(t,\x)\}^2\,dF_{\X,N}(\x)\,dW(t)\to0$.
    For KS,
    $\sup_{t\in\mathcal T,\,\x\in\R^p}|\mu_N(t,\x)|\to0$.
\end{assumption}

Together with Proposition~\ref{prop:fitted-expansions},
Assumption~\ref{ass:bias} implies that $T_N-T_{N,c}=o_p(1)$
for the corresponding CvM or KS statistic.
We next discuss sufficient conditions on $u$ and $N$ for Assumption~\ref{ass:bias}.

By Proposition~\ref{prop:projected-bias-expansion}, Assumption~\ref{ass:bias} holds for both tests if
\begin{equation}
    N^{1/2}
    \sup_{t\in\mathcal T,\,\x\in\R^p}
    \left|\pn\left\{h_{t,\x}^{\perp}(\X;\tht_0)a_{u,t}(\X)\right\}\right|\to0
    \label{eq:general-leading-bias-condition}
\end{equation}
and
\begin{equation}
    N^{1/2}(u^{-2\underline\beta}+u^{-\underline\beta-\underline\eta})\to0.
    \label{eq:hall-wp-threshold}
\end{equation}
These conditions ensure that the leading term and the remainder in $\mu_N(t,\x)$ in \eqref{eq:wp-leading-mean} both converge to zero uniformly over $t\in\mathcal T$ and $\x\in\R^p$.

Under the general Hall class in Assumptions~\ref{ass:hall},
the bound established in Section~\ref{sec:mechanism} gives
\begin{equation}
    \notag 
    \sup_{t\in\mathcal T,\,\x\in\R^p}\left|\pn\left\{h_{t,\x}^{\perp}(\X;\tht_0)a_{u,t}(\X)\right\}\right|=O(u^{-\underline\beta}).
\end{equation}
Consequently, $N^{1/2}u^{-\underline\beta}\to0$, or equivalently $u/N^{1/(2\underline\beta)}\to\infty$, is sufficient for \eqref{eq:general-leading-bias-condition} and
\eqref{eq:hall-wp-threshold}.

If \eqref{eq:regional-bias-rate} holds with
$\kappa>\underline\beta$, a sufficient condition for
Assumption~\ref{ass:bias} for both tests is
$N^{1/2}u^{-q}\to0$ with 
$q=\min\{\kappa,2\underline\beta,\underline\beta+\underline\eta\}$. 
Since $q>\underline\beta$, this condition allows $u$ to increase more slowly relative to $N$ than the condition
$N^{1/2}u^{-\underline\beta}\to0$.

Under \eqref{eq:hall-alignment}, the leading term in $\mu_N(t,\x)$ vanishes, so \eqref{eq:general-leading-bias-condition} holds automatically and \eqref{eq:hall-wp-threshold} alone is sufficient.
The latter is equivalent to $N^{1/2}u^{-2\underline\beta}\to0$ when
$\underline\eta\geq\underline\beta$, and to $N^{1/2}u^{-(\underline\beta+\underline\eta)}\to0$ when
$0<\underline\eta<\underline\beta$.
Both exponents exceed $\underline\beta$, so these conditions can permit a lower threshold than the general requirement
$N^{1/2}u^{-\underline\beta}\to0$, retaining more observations while keeping the Pareto approximation bias after projection asymptotically negligible.

The following theorem establishes the asymptotic size of the CvM and KS tests.

\begin{theorem}
\label{thm:bootstrap}
    For either the CvM or KS statistic, suppose that $H_0$ and
    Assumptions~\ref{ass:hall}, \ref{ass:regular},
    \ref{ass:variation}, and~\ref{ass:bias} hold. 
    Then
    \begin{equation}
        \sup_{z\in\R}
        \left|
        \Prob(T_N\leq z\mid\mathcal C_N)
        -\Prob(T_{N,G}\leq z\mid\mathcal C_N)
        \right|\to0.
        \label{eq:wp-triangular-gaussian}
    \end{equation}
    For $\alpha\in(0,1)$, let $c_{N,1-\alpha}^*$ denote the
    $(1-\alpha)$ quantile of the conditional distribution of $T_N^*$ given the data. 
    Then $\Prob(T_N>c_{N,1-\alpha}^*\mid\mathcal C_N)\to\alpha$.
\end{theorem}

\subsection{Asymptotic size when covariates concentrate}
\label{sec:scope}

Assumption~\ref{ass:regular}(ii) requires the smallest eigenvalue of the information matrix $\Sigma_{\X,0}$ to be bounded away from zero.
As the common threshold $u$ increases, the covariates of observations with $Y_i>u$ may concentrate near values where $\gamma_0(\x)$ is largest.
In log-linear models, this concentration can cause $\lambda_{\min}(\Sigma_{\X,0})$ to approach zero, so Assumption~\ref{ass:regular}(ii) can fail \citep{yangrank2026}.
Under additional conditions, \citet{yangrank2026} establishes asymptotic normality of the parameter estimator under a matrix normalization.
We establish asymptotic size $\alpha$ for the projected CvM and KS bootstrap tests while allowing $\lambda_{\min}(\Sigma_{\X,0})\to0$.
This requires control of the projected test process, the Pareto approximation bias after projection, and the multiplier approximation.

Consider the log-linear specification
\begin{equation}
    \gamma(\x;\tht) = \exp[\{\boldsymbol f(\x)\}\trans\tht],
    \label{eq:scope-loglinear}
\end{equation}
where $\boldsymbol f(\x)\in\R^d$ is a specified vector of regressors and $\sup_{\x\in\mathcal X}\|\boldsymbol f(\x)\|<\infty$.
In this model, $\vv(\x;\tht)=\boldsymbol f(\x)$ and 
$\Sigma_{\X,0} = \pn[\boldsymbol f(\X)\{\boldsymbol f(\X)\}\trans]$. 
When $\Sigma_{\X,0}$ is positive definite, define
\begin{equation}
    \notag
    \lambda_N = \lambda_{\min}(\Sigma_{\X,0})
    \qquad \text{and} \qquad
    \widetilde{\boldsymbol f}_i = \Sigma_{\X,0}^{-1/2}\boldsymbol f(\X_i).
\end{equation}
The normalized regressors satisfy
$N^{-1}\sum_{i=1}^N\widetilde{\boldsymbol f}_i\widetilde{\boldsymbol f}_i\trans=\boldsymbol{\mathsf I}_d$. 
Replacing $\dd_t(\X_i;\tht)$ by $\Sigma_{\X,0}^{-1/2}\dd_t(\X_i;\tht)$ preserves the column space of the derivative matrix and therefore leaves the projection weights in \eqref{eq:weightprojection} unchanged.
We replace Assumption~\ref{ass:regular}(ii) by the following condition.

\begin{assumption}
\label{ass:normalized-design}
    The matrix $\Sigma_{\X,0}$ is positive definite and
    $N^{-1}\sum_{i=1}^N\|\widetilde{\boldsymbol f}_i\|^4=O(1)$. 
    In addition,
    \begin{equation}
        \lambda_N^{-1/2}
        (u^{-\underline\beta}+N^{-1/2})\to0.
        \label{eq:weak-design-rate}
    \end{equation}
\end{assumption}

Assumption~\ref{ass:normalized-design} allows $\lambda_N$ to tend to zero but does not require this convergence.
The normalization gives $N^{-1}\sum_{i=1}^N\|\widetilde{\boldsymbol f}_i\|^2=d$, whereas the fourth-moment bound additionally controls large values of the normalized regressors.
Under Assumption~\ref{ass:regular}(ii), boundedness of $\boldsymbol f$ gives
$N^{-1}\sum_{i=1}^N\|\widetilde{\boldsymbol f}_i\|^4 \leq c_0^{-2}\sup_{\x\in\mathcal X}\|\boldsymbol f(\x)\|^4 < \infty$, 
and \eqref{eq:weak-design-rate} follows from $u,N\to\infty$.
Thus Assumption~\ref{ass:regular}(ii) implies Assumption~\ref{ass:normalized-design} in the bounded log-linear model.

When $\lambda_N\to0$, condition~\eqref{eq:weak-design-rate} requires both $N\lambda_N\to\infty$ and $u^{-\underline\beta}\lambda_N^{-1/2}\to0$.
The two requirements control parameter estimation error and the parameter displacement caused by the Pareto approximation, respectively.
Assumption~\ref{ass:variation} remains a separate condition on the conditional variance of the projected test process. 
Section~A.4 of the Supplementary Material gives sufficient conditions for Assumptions~\ref{ass:normalized-design} and~\ref{ass:variation}.
Under some regularity conditions, $\lambda_N$ decreases at rate $(\log u)^{-2}$, and Assumptions~\ref{ass:normalized-design} and~\ref{ass:variation} hold when $N/(\log u)^2\to\infty$. 

\begin{theorem}
\label{thm:degenerate-design}
    Suppose that $H_0$, Assumptions~\ref{ass:hall},
    \ref{ass:regular}(i), and \ref{ass:normalized-design},
    and the log-linear specification \eqref{eq:scope-loglinear} hold.
    (i) For all sufficiently large $u$ and $N$, $\E\{Q_N(\tht)\mid\mathcal C_N\}$ has a unique local minimizer $\tht_u$ in a sufficiently small neighborhood of $\tht_0$.
    With probability tending to one, $Q_N(\tht)$ has a unique local minimizer $\widehat\tht$ in the same neighborhood.
    These minimizers satisfy
    \begin{equation}
        \|\Sigma_{\X,0}^{1/2}(\tht_u-\tht_0)\|=O(u^{-\underline\beta})
        \qquad \text{~~and~~} \qquad
        \|\Sigma_{\X,0}^{1/2}(\widehat\tht-\tht_u)\|=\Op(N^{-1/2}).
        \label{eq:weak-design-estimation}
    \end{equation}
    (ii) The expansions \eqref{eq:wp-fitted-expansion} and \eqref{eq:wp-leading-mean} hold uniformly over $t\in\mathcal T$ and $\x\in\R^p$.
    (iii) For either the CvM or KS statistic based on
    $U_N(t,\x)-\mu_N(t,\x)$, the Gaussian and multiplier
    approximations in Proposition~\ref{prop:bounded-drift}
    hold under Assumption~\ref{ass:variation}. 
    If Assumption~\ref{ass:bias} also holds, the conclusions of
    Theorem~\ref{thm:bootstrap} apply. 
    In particular, the test using the conditional $(1-\alpha)$ multiplier quantile has asymptotic size $\alpha$.
\end{theorem}

Equation~\eqref{eq:weak-design-estimation} implies the parameter error bounds
$\|\tht_u-\tht_0\|=O(u^{-\underline\beta}\lambda_N^{-1/2})$
and
$\|\widehat\tht-\tht_u\|=\Op\{(N\lambda_N)^{-1/2}\}$. 
Condition~\eqref{eq:weak-design-rate} therefore ensures that
$\tht_u\to\tht_0$ and
$\widehat\tht\overset{p}{\to}\tht_0$,
even when $\lambda_N\to0$. 
By part~(ii) of Theorem~\ref{thm:degenerate-design}, conditions \eqref{eq:general-leading-bias-condition} and
\eqref{eq:hall-wp-threshold} remain sufficient for
Assumption~\ref{ass:bias}.
Under the assumptions of Theorem~\ref{thm:degenerate-design}, $N^{1/2}u^{-\underline\beta}\to0$ is sufficient for Assumption~\ref{ass:bias} for both statistics.
Under \eqref{eq:hall-alignment}, Assumption~\ref{ass:bias} holds under the weaker sufficient condition $N^{1/2}(u^{-2\underline\beta}+u^{-\underline\beta-\underline\eta})\to0$. 
Together with Assumption~\ref{ass:variation}, either condition gives asymptotic size $\alpha$ for both tests by Theorem~\ref{thm:degenerate-design}(iii).

\subsection{Asymptotic power}
\label{sec:power}

In this section, we establish consistency of the proposed projected tests against fixed alternatives and characterize their asymptotic power against local alternatives.

Under a fixed alternative, the true conditional tail index function
$\gamma_0(\x)$ satisfies 
$\Prob(\gamma_0(\X)\neq\gamma(\X;\tht))>0$
for every $\tht\in\Theta$. 
We study power conditional on $\mathcal C_N$. %
As in Section~\ref{sec:estimation}, we write $\tht_u$ for a local
minimizer of $\E\{Q_N(\tht)\mid\mathcal C_N\}$ under the alternative
and $\widehat\tht$ for the local minimizer of $Q_N(\tht)$ near $\tht_u$.
The quantity $\mu_N(t,\x)$ retains its definition in Section~\ref{sec:mechanism} as the conditional mean of the projected test process with the parameter fixed at $\tht_u$. 
Under the fixed alternative, $\mu_N(t,\x)$ can reflect both misspecification of the tail index model and error in the Pareto approximation.

\begin{assumption}
\label{ass:fixed-power}
    (i) For all sufficiently large $u$ and $N$, $\E\{Q_N(\tht)\mid\mathcal C_N\}$ has a local minimizer $\tht_u$.
    There exists a constant $\rho>0$, independent of $u$ and $N$, such that the closed ball of radius $\rho$ centered at $\tht_u$ is contained in the interior of $\Theta$. 
    The boundedness and smoothness conditions in Assumption~\ref{ass:regular}(i) hold uniformly over $\x\in\mathcal X$ and $\tht$ in this ball, with bounds and Lipschitz constants independent of $u$ and $N$. 
    (ii) There exists a constant $c_2>0$ such that, for all sufficiently
        large $u$ and $N$,
        $\lambda_{\min}\left(\left.\partial_{\tht}^{2}\E\{Q_N(\tht)\mid\mathcal  C_N\}\right|_{\tht=\tht_u}\right)\geq c_2 $ and
        $\lambda_{\min}\!\left(\pn[\vv(\X;\tht_u)\{\vv(\X;\tht_u)\}\trans]\right)\geq c_2$.
    (iii) For CvM, 
        $\int_{\mathcal T}\int_{\R^p}\{\mu_N(t,\x)\}^2\,dF_{\X,N}(\x)\,dW(t)\to\infty$.
        For KS,
        $\sup_{t\in\mathcal T,\,\x\in\R^p}|\mu_N(t,\x)|\to\infty$.
\end{assumption}

Together with Assumption~\ref{ass:hall}, parts~(i) and~(ii) of Assumption~\ref{ass:fixed-power} give $\widehat\tht-\tht_u=\Op(N^{-1/2})$ and a uniformly $\Op(1)$ difference between $U_N(t,\x)$ and $\mu_N(t,\x)$.
The conditional multiplier critical values also remain bounded in probability.
Part~(iii) requires the conditional mean to diverge in the integral or supremum used by the corresponding statistic.

\begin{theorem}
\label{prop:wp-fixed-power}
    Suppose that Assumptions~\ref{ass:hall} and~\ref{ass:fixed-power} hold under a fixed alternative. 
    Then, as $u,N\to\infty$,

    \noindent (i) $T_N^{\mathrm{CvM}}\Big/\displaystyle\int_{\mathcal T}\int_{\R^p} \{\mu_N(t,\x)\}^2\,dF_{\X,N}(\x)\,dW(t) \overset{p}{\to} 1$ and 
    $T_N^{\mathrm{KS}}\Big/\displaystyle\sup_{t\in\mathcal T,\,\x\in\R^p}|\mu_N(t,\x)|\overset{p}{\to} 1$.

    \noindent (ii) For every fixed $\alpha\in(0,1)$, the CvM and KS tests are consistent in the sense that
    $\Prob\bigl(T_N^{\mathrm{CvM}}>c_{N,1-\alpha}^{*,\mathrm{CvM}}\bigr)\to1$ and
    $\Prob\bigl(T_N^{\mathrm{KS}}>c_{N,1-\alpha}^{*,\mathrm{KS}}\bigr)\to1$.
\end{theorem}

Section~A.5 of the Supplementary Material distinguishes identification of the tail index function from divergence of the conditional mean $\mu_N(t,\x)$.
At a parameter satisfying the population score equation, projection cannot remove the entire discrepancy between the true and fitted exponential survival functions unless the tail indices agree almost everywhere under the covariate distribution.
When $\nu_u$ converges to a distribution with support $\mathcal X$, Proposition~A.7 gives sufficient conditions for $N^{-1/2}\mu_N(t,\x)$ to have a nonzero limit.
Under those conditions, misspecification on a set of positive limiting covariate probability verifies Assumption~\ref{ass:fixed-power}(iii) and gives consistency of both tests.

When the covariates concentrate, misspecification on regions whose probability under $\nu_u$ tends to zero need not make $\mu_N(t,\x)$ diverge.
For log-linear models, Proposition~A.8 replaces the fixed eigenvalue lower bounds in Assumption~\ref{ass:fixed-power}(ii) by positive definiteness of $\Sigma_{\X,0}$, a bounded fourth moment of the normalized regressors, and $N\lambda_N\to\infty$, while requiring $\tht_u$ to remain in a compact subset of the interior of $\Theta$.
The proposition establishes consistency of the proposed tests when $\lambda_N$ may tend to zero, provided the corresponding divergence condition in Assumption~\ref{ass:fixed-power}(iii) holds.

We next consider local alternatives of the form
\begin{equation*}
    \gamma_N(\x)=\gamma_0(\x)\exp\{k_N^{-1/2}g(\x)\},
    \qquad \gamma_0(\x)=\gamma(\x;\tht_0),
\end{equation*}
where $g$ is a fixed bounded function and $k_N\to\infty$ is deterministic satisfying $N/k_N\overset{p}{\to}1$ as $N\to\infty$.
For power analysis under the local alternative, we condition on $N$ alone.
Conditional on $N$, the $N$ pairs $(Z_i,\X_i)$ are independent draws from the distribution of $(\log(Y/u),\X)$ given $Y>u$. 
As in Section~\ref{sec:estimation}, we take $\tht_u$ and $\widehat\tht$ to be the local minimizers of $\E\{Q_N(\tht)\mid\mathcal C_N\}$ and $Q_N(\tht)$ near $\tht_0$.

\begin{assumption}
\label{ass:local-power}
    (i) 
    The covariate space $\mathcal X\subset\R^p$ is compact and closed under coordinatewise maxima: $\x\vee\boldsymbol y\in\mathcal X$ for all $\x,\boldsymbol y\in\mathcal X$, where $\x\vee\boldsymbol y=(\max\{x_j,y_j\})_{j=1}^{p}$.
    As $N\to\infty$, $\nu_u$ converges weakly to a probability
    measure $\nu$ with continuous marginal distributions and whose support is $\mathcal X$.
    (ii) Assumption~\ref{ass:regular}(i) holds, and $\gamma(\x;\tht)$ and its first two derivatives with respect to $\tht$ are jointly continuous in $(\x,\tht)$ for $\x\in\mathcal X$ and $\tht$ in a neighborhood of $\tht_0$.
    The matrix $\int_{\mathcal X}\vv(\boldsymbol z;\tht_0)\{\vv(\boldsymbol z;\tht_0)\}\trans\,d\nu(\boldsymbol z)$ is positive definite, and $g$ is continuous $\nu$-almost everywhere. 
    (iii) Assumption~\ref{ass:hall} holds with $\gamma_N$ and $b_N$ in place of $\gamma_0$ and $b$, respectively, and with constants and bounds independent of $N$. 
    In \eqref{eq:hall-second-order}, $\beta$ and $\eta$ are fixed, while $c$ and $R$ may depend on $N$.
    Conditions \eqref{eq:general-leading-bias-condition} and \eqref{eq:hall-wp-threshold} hold in probability, with $a_{u,t}$ computed using $\gamma_0$ and $b_N$.
\end{assumption}

Under Assumption~\ref{ass:local-power}(i), $F_{\X,N}$ converges weakly to $\nu$ in probability, so the limiting CvM statistic integrates over $\nu$. 
The support condition in part~(i), together with part~(ii), ensures that the Gaussian counterparts of the CvM and KS statistics are nondegenerate.
The rate conditions in part~(iii) ensure that the contribution of the Pareto approximation error to $\mu_N(t,\x)$ is $\op(1)$ uniformly in $(t,\x)$, even when the second-order coefficient $b_N$ under the local alternative does not converge.

Let $h_{t,\x}^{\perp,\nu}(\cdot;\tht_0)$ denote the population least-squares residual from regressing $h_{\x}(\cdot)$ on the components of $\dd_t(\cdot;\tht_0)$, with the mean squared error evaluated under $\nu$.
Equivalently, replace every empirical average in \eqref{eq:weightprojection} by integration with respect to $\nu$.
Under the local alternative, the limiting mean shift of the projected process is
\begin{equation*}
    \mu_g(t,\x) = \int_{\mathcal X}h_{t,\x}^{\perp,\nu}(\boldsymbol z;\tht_0)\frac{t}{\gamma_0(\boldsymbol z)}S_t(\boldsymbol z;\tht_0)g(\boldsymbol z)\,d\nu(\boldsymbol z).
\end{equation*}
Let $G=\{G(t,\x):(t,\x)\in\mathcal T\times\mathcal X\}$ be a
mean-zero Gaussian process with continuous sample paths and
covariance function
\begin{equation*}
    \Gamma(t,\x;s,\boldsymbol y) = \int_{\mathcal X} h_{t,\x}^{\perp,\nu}(\boldsymbol z;\tht_0) h_{s,\boldsymbol y}^{\perp,\nu}(\boldsymbol z;\tht_0)\{S_{t\vee s}(\boldsymbol z;\tht_0) -S_t(\boldsymbol z;\tht_0)S_s(\boldsymbol z;\tht_0)\} \,d\nu(\boldsymbol z).
\end{equation*}
Let $q_{1-\alpha}^{\mathrm{CvM}}$ and $q_{1-\alpha}^{\mathrm{KS}}$
denote the respective $(1-\alpha)$ quantiles of $\int_{\mathcal T}\int_{\mathcal X}G(t,\x)^2\,d\nu(\x)\,dW(t)$ and $\sup_{t\in\mathcal T,\,\x\in\mathcal X}|G(t,\x)|$.

\begin{theorem}
\label{prop:wp-local-power}
    Suppose that Assumption~\ref{ass:local-power} holds under the local alternatives $\gamma_N(\x)$.
    As $u\to\infty$ and $N/k_N\xrightarrow{p}1$, the following conclusions hold.

    \noindent (i) 
    $(T_N^{\mathrm{CvM}},T_N^{\mathrm{KS}})
        \overset{d}{\to} \Bigl(\int_{\mathcal T}\int_{\mathcal X}\{G(t,\x)+\mu_g(t,\x)\}^2\,d\nu(\x)\,dW(t), \sup_{t\in\mathcal T,\,\x\in\mathcal X} |G(t,\x)+\mu_g(t,\x)|\Bigr)$.

    \noindent (ii) For every fixed $\alpha\in(0,1)$,
    \begin{equation}
        \begin{aligned}
            &\Prob\left(T_N^{\mathrm{CvM}}>c_{N,1-\alpha}^{*,\mathrm{CvM}}\right) \to \Prob\left(
            \int_{\mathcal T}\int_{\mathcal X}
            \{G(t,\x)+\mu_g(t,\x)\}^2\,d\nu(\x)\,dW(t)
            >q_{1-\alpha}^{\mathrm{CvM}}
            \right),\\
            &\Prob\left(T_N^{\mathrm{KS}}>c_{N,1-\alpha}^{*,\mathrm{KS}}\right) \to \Prob\left(
            \sup_{t\in\mathcal T,\,\x\in\mathcal X}
            |G(t,\x)+\mu_g(t,\x)|
            >q_{1-\alpha}^{\mathrm{KS}}
            \right).
        \end{aligned}
        \label{eq:wp-local-power}
    \end{equation}

    \noindent (iii) If there exists $\boldsymbol a\in\R^d$ such
    that $g(\x)=\{\vv(\x;\tht_0)\}\trans\boldsymbol a$ for
    $\nu$-almost every $\x$, both limits in
    \eqref{eq:wp-local-power} equal $\alpha$. Otherwise, both
    limits lie strictly between $\alpha$ and one.
\end{theorem}

Theorem~\ref{prop:wp-local-power} shows that the proposed tests have nontrivial asymptotic power against local alternatives converging to the null at rate $N^{-1/2}$ whenever $g$
lies outside the tangent space of the null model, which consists
of the linear combinations of the components of $\vv(\cdot;\tht_0)$.
If $g(\cdot)=\{\vv(\cdot;\tht_0)\}\trans\boldsymbol a$ for some $\boldsymbol a\in\R^d$, the local alternative agrees to first order with the null model at
$\tht_0+k_N^{-1/2}\boldsymbol a$. 
Projection then gives $\mu_g(t,\x)=0$, and the proposed tests have limiting rejection probability $\alpha$.

\section{Threshold selection}
\label{sec:threshold}

We propose a novel threshold selection rule based on the estimated Pareto approximation bias after projection.
The objective is to retain more observations for testing while limiting the estimated bias relative to the bootstrap critical values, with separate criteria for the CvM and KS tests.

We use the full sample of $n$ observations for threshold selection.
Let $\widehat u_{\mathrm{ref}}$ be a reference threshold obtained, for example, by the method of \citet{wangtsai2009}.
Consider a finite, nonempty set of candidate thresholds $\mathcal W_n$ with $\max_{u \in \mathcal W_n} u=\widehat u_{\mathrm{ref}}$.
For example, $\mathcal W_n$ can be a grid in 
$[\max\{u_{\mathrm{floor}},c_3\widehat u_{\mathrm{ref}}\}, \widehat u_{\mathrm{ref}}]$, 
including $\widehat u_{\mathrm{ref}}$, where $c_3\in(0,1)$ and $0<u_{\mathrm{floor}}\leq\widehat u_{\mathrm{ref}}$.
The lower bound $u_{\mathrm{floor}}$ is chosen to keep the number of observations above each candidate threshold within a prescribed upper limit. 
For $u\in\mathcal W_n$, let $K_u=\sum_{i=1}^n\ind\{Y_i>u\}$ and write $\x_1,\ldots,\x_{K_u}$ for the covariates of the observations with $Y_i>u$. 
Thus the sample size used for testing at threshold $u$ is $N=K_u$.

The conditional mean of the projected test process can reflect both Pareto approximation error and misspecification of the conditional tail index model.
To assess the Pareto approximation bias, we compare the conditional means of the test process under the true conditional distribution of the log excess and an exponential distribution with mean $\gamma_0(\x)$.
The exponential distribution corresponds to an exact Pareto tail with the same conditional tail index. 
Subtracting the conditional mean under the exponential distribution avoids treating model misspecification as Pareto approximation bias. 

Let $\widehat F(y\mid\x)$ estimate the conditional distribution function of $Y$ given $\X=\x$, and let $\widehat\gamma(\x)$ be an estimate of $\gamma_0(\x)$ obtained without imposing $H_0$.
For each candidate threshold $u$, the estimated conditional survival function and conditional mean of $Z=\log(Y/u)$ given $Y>u$ and $\X=\x$ are
\begin{equation}
    \notag
    \widehat p_u^{\mathrm{exc}}(t,\x)
    =
    \frac{1-\widehat F(ue^t\mid\x)}
         {1-\widehat F(u\mid\x)},
    \qquad
    \widehat\mu_u^{\mathrm{exc}}(\x)
    =
    \int_0^\infty\widehat p_u^{\mathrm{exc}}(t,\x)\,dt.
\end{equation}
For the exponential comparison, define
$\widehat p_u^{\mathrm{exp}}(t,\x) = \exp\{-t/\widehat\gamma(\x)\}$ and
$\widehat\mu_u^{\mathrm{exp}}(\x) = \widehat\gamma(\x)$.

For each $r\in\{\mathrm{exc},\mathrm{exp}\}$, obtain $\widehat\tht_u^r$ by minimizing the negative log quasi-likelihood function  
$K_u^{-1}\sum_{i=1}^{K_u}\bigl\{\log\gamma(\x_i;\tht)+\widehat\mu_u^r(\x_i)/\gamma(\x_i;\tht)\bigr\}$ 
over $\tht\in\Theta$. 
Compute the projection weights $h_{t,\x}^{\perp}(\x_i;\widehat\tht_u^r)$ from \eqref{eq:weightprojection}, with empirical averages taken over $\x_1,\ldots,\x_{K_u}$.

For $(t,\x)\in\mathcal T\times\R^p$, define the estimated Pareto approximation bias after projection by
\begin{equation}
    \notag
    \begin{aligned}
        \widehat B_u(t,\x)
        ={}&
        K_u^{-1}\sum_{i=1}^{K_u}
        h_{t,\x}^{\perp}(\x_i;\widehat\tht_u^{\mathrm{exc}})
        \left\{
            \widehat p_u^{\mathrm{exc}}(t,\x_i)
            -
            S_t(\x_i;\widehat\tht_u^{\mathrm{exc}})
        \right\}\\
        &-
        K_u^{-1}\sum_{i=1}^{K_u}
        h_{t,\x}^{\perp}(\x_i;\widehat\tht_u^{\mathrm{exp}})
        \left\{
            \widehat p_u^{\mathrm{exp}}(t,\x_i)
            -
            S_t(\x_i;\widehat\tht_u^{\mathrm{exp}})
        \right\}.
    \end{aligned}
\end{equation}
Both averages use the same covariate values, but each uses the parameter estimates and projection weights obtained under its conditional distribution.

Under $H_0$, replacing $\widehat F(y\mid \x)$ and $\widehat\gamma(\x)$ by $F(y\mid \x)$ and $\gamma_0(\x)$ gives the expected quasi-likelihood fits $\tht_u$ and $\tht_0$, respectively.
The exponential residual mean is then zero because $\exp\{-t/\gamma_0(\x_i)\}=S_t(\x_i;\tht_0)$.
Consequently, the difference of the weighted averages equals $K_u^{-1/2}\mu_N(t,\x)$, with $N=K_u$.
Under an alternative, the exponential weighted residual mean need not vanish and is subtracted to exclude the contribution present under an exact Pareto model. 
The estimated exponential contribution can also be nonzero under $H_0$, because $\widehat\gamma(\x)$ is not restricted to the null model.

Since the test process is $K_u^{1/2}$ times an average of weighted residuals, we use $K_u^{1/2}\widehat B_u(t,\x)$ as the estimate of the Pareto approximation bias on the scale of the test process.
At each candidate threshold, fit the null model to the observed log excesses with $Y_i>u$ and apply the multiplier bootstrap in Section~\ref{sec:bootstrap}. 
Let $\widehat c_C(u)>0$ and $\widehat c_K(u)>0$ denote the conditional $(1-\alpha)$ quantiles of the bootstrap CvM and KS statistics, respectively.

Define the threshold selection criteria by 
\begin{equation}
    \widehat D_C(u) = \frac{K_u^{1/2}\left[
        K_u^{-1}\sum_{a=1}^{K_u}\int_{\mathcal T}
        \{\widehat B_u(t,\x_a)\}^2\,dW(t)
        \right]^{1/2}}
        {\widehat c_C(u)^{1/2}}
        \quad \text{and} \quad
        \widehat D_K(u) = \frac{K_u^{1/2}
        \sup_{t\in\mathcal T,\,\x\in\R^p}|\widehat B_u(t,\x)|}
        {\widehat c_K(u)}.
    \label{ad-eq:criteria}
\end{equation}

For a prescribed tolerance $\eta_n>0$, select
\begin{equation}
    \widehat u_j
    =
    \min\left(
        \left\{
            u\in\mathcal W_n:
            \widehat D_j(u)\leq\eta_n
        \right\}
        \cup
        \{\widehat u_{\mathrm{ref}}\}
    \right),
    \qquad j\in\{C,K\}.
    \label{ad-eq:separate-selectors}
\end{equation}
For each test, \eqref{ad-eq:separate-selectors} selects the lowest candidate satisfying $\widehat D_j(u)\leq\eta_n$, thereby retaining the largest number of observations among the qualifying candidates.
If no candidate satisfies the criterion, the reference threshold is retained.
The CvM and KS thresholds are selected separately and may differ. 
Our simulation studies suggest that $\eta_n=1/\log n$ is a robust choice, with empirical rejection frequencies generally close to the nominal level even for models that do not satisfy \eqref{eq:hall-alignment}. 
Section~A.6 of the Supplementary Material describes the estimators $\widehat F(y\mid \x)$ and $\widehat\gamma(\x)$, the reference threshold, and the implementation of the threshold selection rule.

\section{Simulation studies}
\label{sec:simulations}

In this section, we conduct a series of Monte Carlo experiments to study the finite sample performance of our proposed projected CvM and KS tests. 
Under the null, we compare the size of the projected and unprojected CvM and KS tests. 
Under the alternatives, we compare the power of the projected tests using two rules for threshold selection, one proposed by \citet{wangtsai2009} and the other introduced in Section~\ref{sec:threshold}.
Section~A.1 of the Supplementary Material gives implementation details for the projected and unprojected tests. 

\subsection{Simulation design}

We test whether the log conditional tail index is linear in the covariates, with an intercept and all $p$ covariates included in the fitted model. 
We consider sample sizes $n\in\{5{,}000,10{,}000\}$ and covariate dimensions $p\in\{4,8\}$. 
For each observation, we draw $\boldsymbol G=(G_1,\ldots,G_p)\trans\sim N_p(\boldsymbol 0,\Sigma)$ with $\Sigma=\bigl(0.5^{|j-k|}\bigr)_{1\leq j,k\leq p}$ and set $X_j=12^{1/2}\{\Phi(G_j)-1/2\}$ for $j=1,\ldots,p$, where $\Phi(\cdot)$ is the standard normal distribution function. 
Each $X_j$ is thus uniform on $[-3^{1/2},3^{1/2}]$, with mean zero and variance one. 
In the main paper, we consider the conditional tail index function of the form
\begin{equation}
    \log\gamma_0(\X)
    =\log(0.5)+0.25(-X_1+0.7X_2-0.8X_3+0.5X_4)
    +\delta(X_1^2-1).
    \label{eq:new-numerical-index}
\end{equation}
The null hypothesis holds at $\delta=0$. For the power comparisons, we take $\delta\in\{0.1,0.2,0.3,0.4\}$.
Since $\E(X_1^2)=1$, the quadratic term changes the shape of the log tail index function while leaving its unconditional mean unchanged.
For $p=8$, the four additional covariates have zero coefficients.
We also consider another alternative obtained by replacing $\delta(X_1^2-1)$ with $\delta\sin(\pi X_1/3^{1/2})$.
Section~C.4 of the Supplementary Material reports the empirical power of the projected tests under this alternative.

The first three models belong to the Burr XII family with conditional survival function
\begin{equation*}
    \overline F(y\mid\X)
    =\left[1+\left\{y/s(\X)\right\}^{\beta(\X)}
      \right]^{-1/\{\beta(\X)\gamma_0(\X)\}},
    \qquad y>0.
\end{equation*}
\begin{enumerate}[label={},leftmargin=0pt,itemindent=0pt,labelsep=0pt]
    \item \textit{YI.} The scale and shape are constant, with $s(\X)=1$ and $\beta(\X)=2$.
    \item \textit{YI-scale.} We vary the scale with $X_3$ by setting $s(\X)=\exp(-0.60X_3)$ and retain $\beta(\X)=2$.
    \item \textit{YI-shape.} We retain $s(\X)=1$ and vary the shape with $X_1$ by setting $\beta(\X)=2\exp(0.8X_1/3^{1/2})$.
    \item \textit{YII (Student-$t$).} The conditional survival function is
    $\overline F(y\mid\X)=\bigl\{1-y/(y^2+2)^{1/2}\bigr\}^{1/\{2\gamma_0(\X)\}}$, $y>0$.
    \item \textit{YIII (log-logistic).} The conditional distribution is supported on $[1,\infty)$, with
    $\overline F(y\mid\X)=1.3\big/\bigl\{0.3+y^{1/\gamma_0(\X)}\bigr\}$, $y\geq1$.
    \item \textit{YV (Fr\'echet).} We use a Fr\'echet distribution with unit scale and conditional survival function
    $\overline F(y\mid\X)=1-\exp\{-y^{-1/\gamma_0(\X)}\}$, $y>0$.
\end{enumerate}

YI and YI-scale illustrate how the scale of the conditional distribution affects the Pareto approximation bias after projection.
In both models, the second-order coefficient in Assumption~\ref{ass:hall} satisfies $\gamma_0(\X)b(\X)=-s(\X)^2/2$.
In YI, this quantity is constant, so the intercept suffices for \eqref{eq:hall-alignment} to hold under the null. YII satisfies the same condition.
Projection therefore removes the leading Pareto approximation bias in both models.
In YI-scale, $-s(\X)^2/2$ is nonlinear in $X_3$ and lies outside the span of the regressors, so the leading Pareto approximation bias need not vanish after projection.

We consider three threshold rules. FIX uses the empirical $(1-\rho)$ quantile of $Y$ in the full sample.
WT denotes the method of \citet{wangtsai2009}, with the search restricted to thresholds that retain at most a fraction $\rho$ of this sample. 
AD applies the threshold selection rule in Section~\ref{sec:threshold} with tolerance $\eta_n=1/\log n$ on a grid whose largest threshold is the WT reference. 
FIX and WT each use a single threshold for the CvM and KS tests, both with and without projection. 
AD selects thresholds separately for the projected CvM and KS tests, for which its results are reported.
We use $\rho\in\{20\%,25\%,30\%\}$ for size and $\rho=30\%$ for power comparison.
The nominal significance level is $5\%$. 
All results are based on 500 replications.

\subsection{Size}

Table~\ref{tab:new-size} reports null rejection frequencies for $n=10{,}000$ and $p=4$. 
Section~C.1 of the Supplementary Material gives the remaining size results. 
Under FIX, the unprojected CvM and KS tests overreject in all six models, and their rejection frequencies increase with $\rho$. 
WT reduces this distortion for the unprojected tests, although their rejection frequencies remain above the nominal level in several settings. The AD results are reported for the projected tests only.

In YI and YII, condition~\eqref{eq:hall-alignment} holds, so projection removes the leading Pareto approximation bias. 
In both models, the projected tests have rejection frequencies close to $5\%$ even at $\rho=30\%$, the lowest of the three FIX thresholds. 
In YI-scale, condition~\eqref{eq:hall-alignment} fails, and lowering the threshold can lead to overrejection even after projection. At $\rho=30\%$, the projected CvM and KS tests reject in $24.6\%$ and $24.4\%$ of replications under FIX.

The YI-scale results illustrate the importance of data-driven threshold selection when the leading Pareto approximation bias need not vanish after projection. 
WT selects the threshold by assessing the fit of the exponential approximation to the log excesses. 
AD instead uses $\widehat B_u$ to estimate the Pareto approximation bias after projection and applies the test-specific scaling and normalization in \eqref{ad-eq:criteria}. 
Both rules keep the rejection frequencies of the projected tests close to $5\%$ in YI-scale. 
Under AD, they range from $3.8\%$ to $4.8\%$ across all three values of $\rho$. 
We therefore restrict the power comparisons to the projected CvM and KS tests under WT and AD. 

\begin{table}[!htbp]\centering
\caption{Null rejection frequencies (\%) for $n=10{,}000$ and $p=4$.}
\label{tab:new-size}
\begingroup\linespread{1}\fontsize{9.3}{10.8}\selectfont
\setlength{\tabcolsep}{1.6pt}\renewcommand{\arraystretch}{1.00}\setlength{\aboverulesep}{.23ex}\setlength{\belowrulesep}{.32ex}
\begin{tabular*}{\linewidth}{@{\extracolsep{\fill}}llrrrrrrrrrr@{}}
\toprule
 &  & \multicolumn{6}{c}{Projected} & \multicolumn{4}{c}{Unprojected}\\
\cmidrule(lr){3-8}\cmidrule(lr){9-12}
 &  & \multicolumn{2}{c}{FIX} & \multicolumn{2}{c}{WT} & \multicolumn{2}{c}{AD} & \multicolumn{2}{c}{FIX} & \multicolumn{2}{c}{WT}\\
\cmidrule(lr){3-4}\cmidrule(lr){5-6}\cmidrule(lr){7-8}\cmidrule(lr){9-10}\cmidrule(lr){11-12}
Model & $\rho$ & CvM & KS & CvM & KS & CvM & KS & CvM & KS & CvM & KS\\
\midrule
YI & 20\% & 6.2 & 6.0 & 4.6 & 6.0 & 3.2 & 5.8 & 45.8 & 54.4 & 4.4 & 5.4\\
 & 25\% & 4.2 & 5.8 & 4.6 & 5.8 & 3.4 & 5.6 & 75.0 & 83.0 & 5.4 & 5.8\\
 & 30\% & 5.2 & 5.4 & 4.4 & 5.8 & 3.4 & 5.6 & 95.8 & 98.6 & 5.6 & 5.8\\
\cmidrule(l){2-12}
YI-scale & 20\% & 10.4 & 6.8 & 3.4 & 5.2 & 3.8 & 4.8 & 98.6 & 99.2 & 13.2 & 11.8\\
 & 25\% & 15.8 & 14.2 & 3.4 & 5.2 & 3.8 & 4.8 & 100.0 & 100.0 & 13.2 & 11.8\\
 & 30\% & 24.6 & 24.4 & 3.4 & 5.2 & 3.8 & 4.8 & 100.0 & 100.0 & 13.2 & 11.8\\
\cmidrule(l){2-12}
YI-shape & 20\% & 3.8 & 5.4 & 3.2 & 6.6 & 3.6 & 6.0 & 67.0 & 67.4 & 38.2 & 33.2\\
 & 25\% & 5.0 & 7.2 & 3.0 & 6.0 & 3.6 & 6.0 & 88.2 & 90.2 & 48.0 & 44.8\\
 & 30\% & 5.8 & 5.2 & 3.2 & 5.2 & 3.6 & 4.8 & 98.0 & 98.8 & 52.2 & 49.6\\
\cmidrule(l){2-12}
YII & 20\% & 3.8 & 5.2 & 4.4 & 4.6 & 3.8 & 5.6 & 79.2 & 88.4 & 8.0 & 7.0\\
 & 25\% & 4.6 & 5.4 & 4.4 & 4.6 & 3.8 & 5.6 & 98.6 & 99.6 & 8.0 & 7.0\\
 & 30\% & 5.4 & 6.0 & 4.4 & 4.6 & 3.8 & 5.6 & 100.0 & 100.0 & 8.0 & 7.0\\
\cmidrule(l){2-12}
YIII & 20\% & 6.8 & 5.8 & 5.6 & 5.6 & 4.2 & 5.4 & 8.6 & 11.2 & 4.4 & 3.0\\
 & 25\% & 5.2 & 5.8 & 6.2 & 4.4 & 4.2 & 4.6 & 10.0 & 11.0 & 4.4 & 2.6\\
 & 30\% & 5.2 & 4.8 & 6.4 & 5.0 & 4.2 & 5.0 & 12.6 & 15.8 & 4.4 & 2.8\\
\cmidrule(l){2-12}
YV & 20\% & 4.2 & 2.8 & 4.6 & 4.6 & 5.2 & 4.8 & 19.8 & 23.2 & 7.2 & 5.6\\
 & 25\% & 3.8 & 5.2 & 5.4 & 6.2 & 5.2 & 5.0 & 36.0 & 42.8 & 8.8 & 7.6\\
 & 30\% & 5.4 & 5.0 & 5.4 & 6.2 & 5.2 & 4.8 & 57.0 & 65.4 & 9.6 & 8.8\\
\bottomrule\end{tabular*}\endgroup
\end{table}

\subsection{Empirical power}
\label{sec:threshold-illustration}

\begin{figure}[htp]
    \centering
    \includegraphics[width=\linewidth]{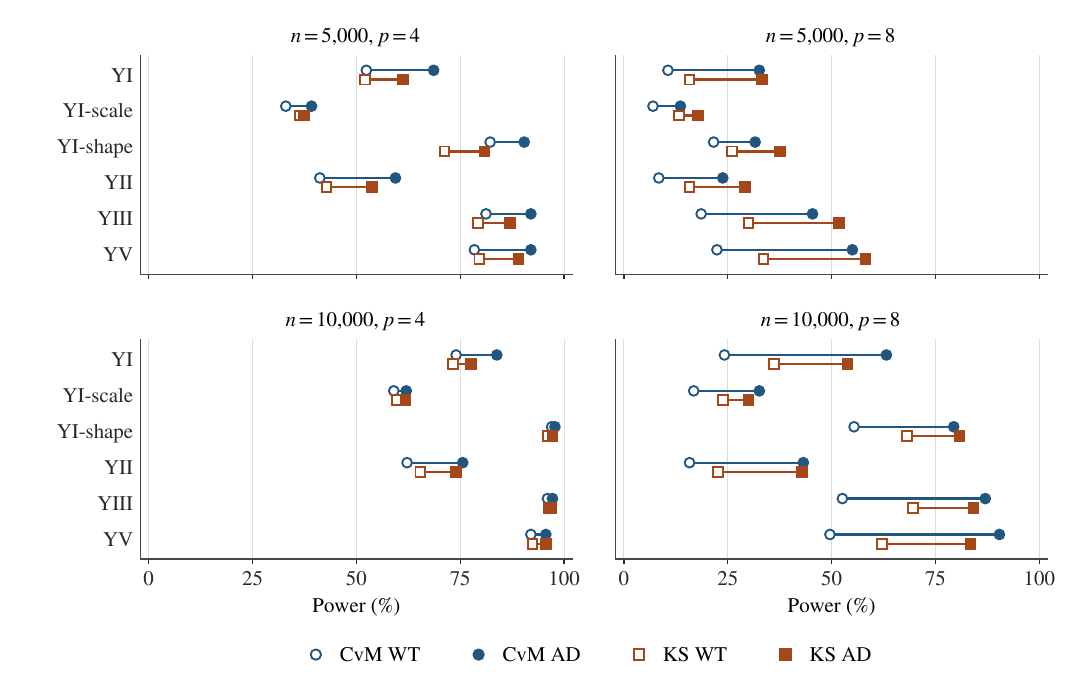}
    \caption{Power of the projected tests for $\delta=0.2$. Open symbols indicate WT and filled symbols indicate AD.}
    \label{fig:power-strength}
\end{figure}

We compare the empirical power of the projected CvM and KS tests under WT and AD with $\rho=30\%$. 
Figure~\ref{fig:power-strength} reports the results for $\delta=0.2$. 
Both tests have higher power under AD than under WT in all six models and for all four combinations of $n$ and $p$. 
The gains are generally largest when $n=10{,}000$ and $p=8$. 
At $n=10{,}000$ and $p=4$, WT already has high power in several models, leaving less room for improvement. 
For both tests, the power gains are generally smaller in YI-scale, where condition~\eqref{eq:hall-alignment} need not hold. 

\begin{figure*}[htp]
    \centering\includegraphics[width=\linewidth]{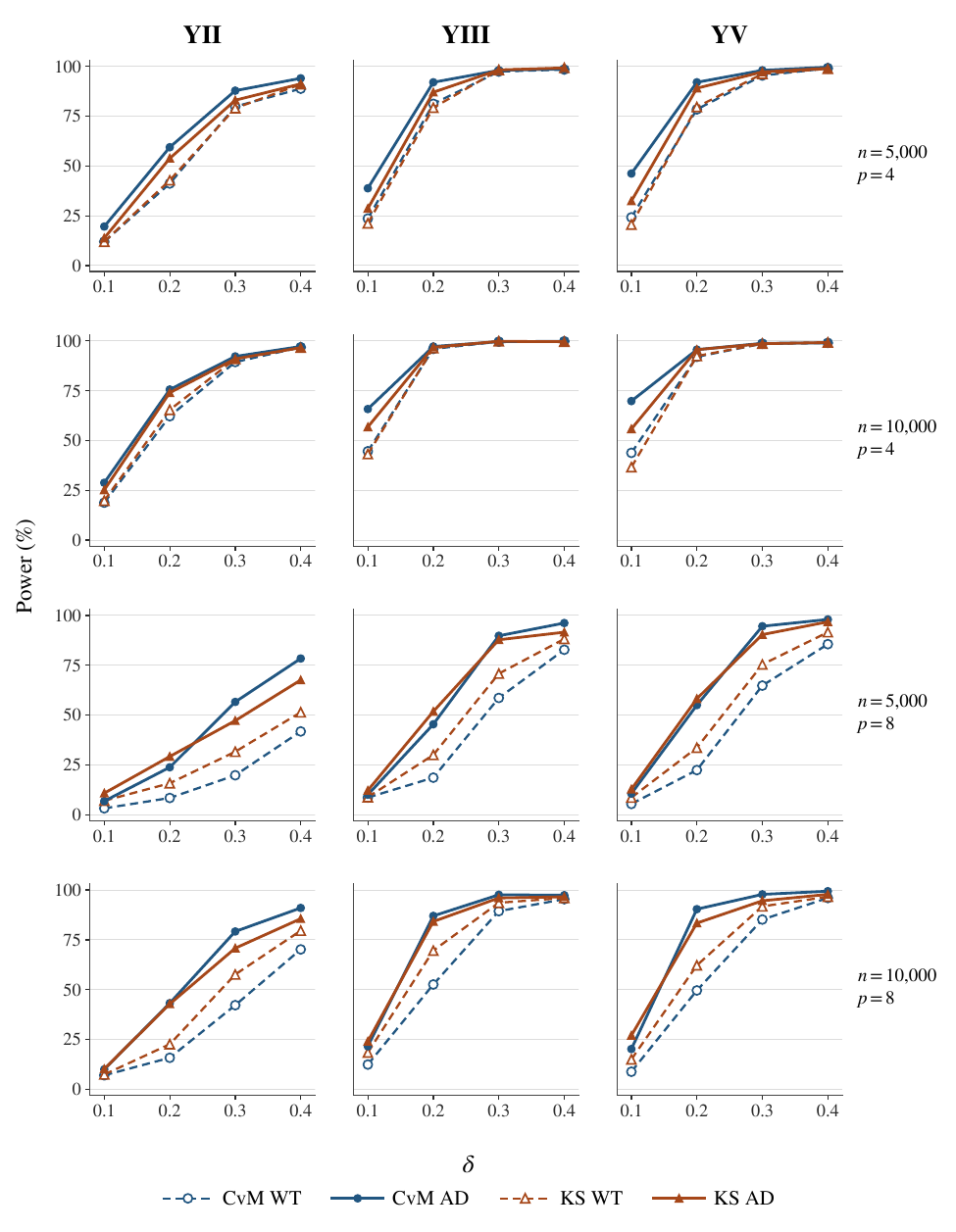}
    \caption{Power of the projected tests in YII, YIII and YV as $\delta$ varies.}
    \label{fig:power-curves-main}
\end{figure*}

Figure~\ref{fig:power-curves-main} shows how power varies with the coefficient $\delta$ of the quadratic term in YII, YIII and YV. Power increases with $\delta$ under both threshold rules. 
For $n=10{,}000$ and $p=4$, the differences between the two rules are small by $\delta=0.3$, when both WT tests have power of at least $95\%$ in all three models.
With eight covariates, AD continues to yield substantial gains at $\delta=0.3$. 
Table~S3 of the Supplementary Material reports the complete power results. 

Table~\ref{tab:main-v11-counts-10000} compares the effective sample sizes under WT and AD for $n=10{,}000$. 
For both projected tests, AD retains more observations on average than WT at every combination of model, $p$ and $\delta$.
In YV with $p=8$ and $\delta=0.2$, WT uses $1{,}823$ observations on average, compared with $3{,}512$ for AD with CvM and $2{,}833$ for AD with KS. 
The CvM and KS rejection frequencies are 49.6\% and 62.2\% under WT and 90.4\% and 83.4\% under AD. 
Section~C.3 of the Supplementary Material reports the corresponding results for $n=5{,}000$.

\begin{table*}[!ht]\centering
\caption{Mean numbers of observations above the selected thresholds for the projected tests ($n=10{,}000$).}
\label{tab:main-v11-counts-10000}
\begingroup\linespread{1}\fontsize{9.5}{11.5}\selectfont
\setlength{\tabcolsep}{2.0pt}\renewcommand{\arraystretch}{1.02}
\begin{tabular*}{\linewidth}{@{\extracolsep{\fill}}llrrrrrrrrrr@{}}
\toprule
 & & \multicolumn{5}{c}{$p=4$} & \multicolumn{5}{c}{$p=8$}\\
\cmidrule(lr){3-7}\cmidrule(lr){8-12}
Model & Rule & $\delta=0$ & $0.1$ & $0.2$ & $0.3$ & $0.4$
             & $\delta=0$ & $0.1$ & $0.2$ & $0.3$ & $0.4$\\
\midrule
YI & WT & 1299 & 1291 & 1258 & 1278 & 1276 & 1296 & 1284 & 1293 & 1270 & 1255\\
 & AD (CvM) & 2596 & 2085 & 1582 & 1430 & 1364 & 3040 & 2952 & 2600 & 2226 & 2017\\
 & AD (KS) & 2438 & 1874 & 1407 & 1307 & 1281 & 2964 & 2555 & 1916 & 1558 & 1423\\
\cmidrule(l){2-12}
YI-scale & WT & 909 & 935 & 959 & 1046 & 1119 & 955 & 910 & 971 & 1056 & 1089\\
 & AD (CvM) & 1124 & 1095 & 1033 & 1073 & 1132 & 1769 & 1624 & 1540 & 1471 & 1387\\
 & AD (KS) & 997 & 1024 & 992 & 1058 & 1125 & 1338 & 1255 & 1185 & 1170 & 1153\\
\cmidrule(l){2-12}
YI-shape & WT & 1946 & 2073 & 2425 & 2322 & 2105 & 1985 & 2144 & 2455 & 2475 & 2132\\
 & AD (CvM) & 2440 & 2707 & 2807 & 2431 & 2213 & 3500 & 3552 & 3655 & 3394 & 2931\\
 & AD (KS) & 2405 & 2788 & 2856 & 2402 & 2134 & 2978 & 3189 & 3271 & 2898 & 2363\\
\cmidrule(l){2-12}
YII & WT & 983 & 1042 & 993 & 1006 & 1040 & 1056 & 1047 & 991 & 1048 & 1069\\
 & AD (CvM) & 2024 & 1716 & 1256 & 1148 & 1120 & 2498 & 2470 & 2004 & 1811 & 1683\\
 & AD (KS) & 1926 & 1547 & 1125 & 1044 & 1048 & 2410 & 2141 & 1513 & 1301 & 1210\\
\cmidrule(l){2-12}
YIII & WT & 2282 & 2261 & 2203 & 2175 & 2203 & 2269 & 2267 & 2207 & 2248 & 2158\\
 & AD (CvM) & 3581 & 3305 & 2823 & 2540 & 2400 & 3778 & 3790 & 3604 & 3535 & 3229\\
 & AD (KS) & 3458 & 2994 & 2474 & 2241 & 2220 & 3707 & 3443 & 2977 & 2716 & 2434\\
\cmidrule(l){2-12}
YV & WT & 1850 & 1847 & 1799 & 1816 & 1769 & 1953 & 1910 & 1823 & 1867 & 1882\\
 & AD (CvM) & 3458 & 3149 & 2495 & 2136 & 1904 & 3782 & 3752 & 3512 & 3352 & 3020\\
 & AD (KS) & 3244 & 2833 & 2168 & 1921 & 1788 & 3691 & 3489 & 2833 & 2515 & 2217\\
\bottomrule
\end{tabular*}\endgroup
\end{table*}

\section{Real data application}
\label{sec:real-data}

We apply the proposed specification tests to first-quarter losses in the \textit{French Commercial Fire Insurance dataset}, available as \texttt{frecomfire} in the \texttt{CASdatasets} R package.
We use the 2,373 commercial fire losses collected by the FFSA in first-quarters of 1982--1996.
The response is the reported loss amount, expressed in thousands of euros at 2007 prices.

\begin{figure}[!htbp]
    \centering
    \includegraphics[width=0.48\textwidth]{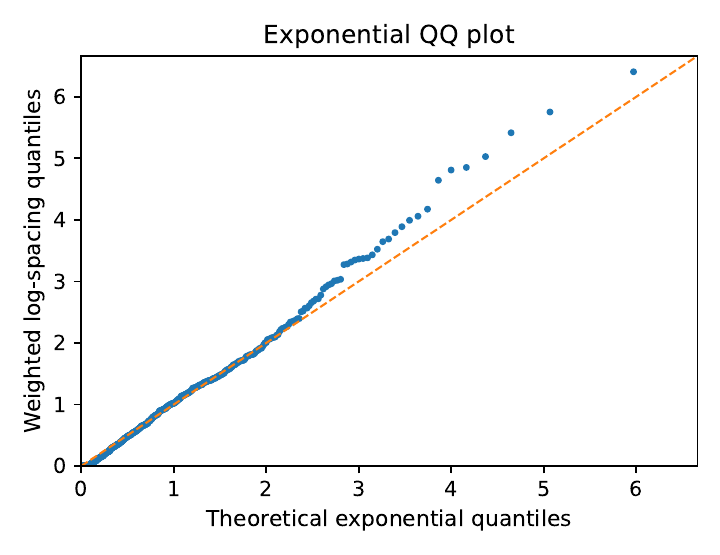}
    \caption{Exponential quantile--quantile plot of the first $\lfloor0.30n\rfloor=711$ weighted log spacings of first-quarter commercial fire losses in millions of euros.
    The dashed line is the identity line.}
    \label{fig:french-fire-tail}
\end{figure}

We first examine the upper tail of the loss distribution.
Let $Y_{(1)}\leq\cdots\leq Y_{(n)}$ denote the ordered losses.
We examine the weighted log spacings $i\log\{Y_{(n-i+1)}/Y_{(n-i)}\}$ for $i=1,\ldots,\lfloor\rho n\rfloor$. 
For a Pareto-type upper tail, the weighted log spacings for small $i$ are approximately independent exponential random variables.
We compare the ordered spacings with the quantiles of an exponential distribution whose mean is estimated by the average of the spacings.
We examine $\rho\in\{0.5,0.4,0.3,0.2\}$.
The quantile--quantile plot for $\rho=0.30$ in Figure~\ref{fig:french-fire-tail} showed an approximately linear relationship, which further confirmed the heavy--tailedness of the loss distribution.
We take $\rho=0.30$, giving $K_{\max}=\lfloor\rho n\rfloor=711$.

Let $X_1$ denote the occurrence year, and let $X_2$ denote the occurrence month.
January, February, and March are coded as 1, 2, and 3, respectively.
Both covariates are standardized using their full-sample means and sample standard deviations.
We test the adequacy of two parametric specifications for the conditional tail index of first-quarter commercial fire losses: a log-linear model in $X_1$ and $X_2$, and an expanded model that also includes their interaction and a quadratic term in $X_1$.
\begin{equation}
    \notag
    \begin{aligned}
        \text{Log-linear model}\qquad
        &\log\gamma(\x;\tht)
        =\beta_0+\beta_1x_1+\beta_2x_2,\\
        \text{Expanded model}\qquad
        &\log\gamma(\x;\tht)
        =\beta_0+\beta_1x_1+\beta_2x_2
         +\beta_3x_1x_2+\beta_4x_1^2.
    \end{aligned}
\end{equation}

For each specification, we compute the projected CvM and KS statistics under the WT and AD threshold rules.
In this application, $K_{\min}=\max\{20d,\lceil n^{1/2}\rceil\}$, giving lower count limits of 60 and 100 for the log-linear and expanded models, respectively.
WT searches every integer $K$ between $K_{\min}$ and $K_{\max}$.
AD evaluates every distinct order-statistic threshold in $[\max\{u_{\mathrm{floor}},\widehat u_{\mathrm{ref}}/2\},\widehat u_{\mathrm{ref}}]$, together with the lower endpoint, where $u_{\mathrm{floor}}=Y_{(n-K_{\max})}$.
We use the tolerance $\eta_n=1/\log n$ and the remaining implementation details in Section~A.6 of the Supplementary Material.
Threshold selection uses 999 independent standard normal multiplier draws, and the reported bootstrap $p$-values are obtained from a separate set of 999 independent standard normal multiplier draws.

Table~\ref{tab:french-fire-current} gives the selected thresholds, the numbers of observations with $Y_i>u$, and the corresponding $p$-values.
For the log-linear model, neither test rejects at the 5\% level under WT (CvM: $p=0.063$; KS: $p=0.152$), whereas the CvM test rejects under AD ($p=0.036$) and the KS test does not ($p=0.101$).
For both tests, AD selects $u=1.943$ million euros, compared with the WT threshold of $2.079$ million euros.
The lower threshold increases the number of observations used for testing from 450 to 486. 
For the expanded model, neither test rejects under WT %
or AD. %
AD selects $u=1.870$ million euros for the expanded CvM test and $u=1.943$ million euros for the expanded KS test, retaining 509 and 486 observations, respectively.
These results indicate that the two tests detect no evidence against the expanded model. 

\begin{table}[!htbp]
    \centering
    \caption{Projected specification tests for first-quarter commercial fire losses.
    Thresholds $u$ are in millions of euros.}
    \label{tab:french-fire-current}
    \small
    \begin{tabular}{@{}lllrrr@{}}
        \toprule
        Model & Test & Rule & $u$ & $K_u$ & $p$-value\\
        \midrule
        Log-linear & CvM & WT & 2.079 & 450 & 0.063\\
         &  & AD & 1.943 & 486 & \textbf{0.036}\\
        \addlinespace
         & KS & WT & 2.079 & 450 & 0.152\\
         &  & AD & 1.943 & 486 & 0.101\\
        \midrule
        Expanded & CvM & WT & 2.079 & 450 & 0.408\\
         &  & AD & 1.870 & 509 & 0.561\\
        \addlinespace
         & KS & WT & 2.079 & 450 & 0.135\\
         &  & AD & 1.943 & 486 & 0.082\\
        \bottomrule
    \end{tabular}
\end{table}

\section{Conclusion}
\label{sec:conclusion}

We have developed projected CvM and KS specification tests for
parametric conditional tail index models, with multiplier bootstrap
critical values.
Projection removes the first-order effect of parameter estimation
and, within the conditional Hall class, eliminates the leading
Pareto approximation bias under condition~\eqref{eq:hall-alignment}.
This cancellation permits lower thresholds while preserving
asymptotically correct size.
We establish power properties against fixed and local alternatives
and extend the size results to log-linear models with covariate
concentration under the stated design conditions.
The proposed threshold selection rule uses the estimated Pareto
approximation bias after projection to determine whether a reference threshold can be lowered.
A natural future direction is to extend the tests to thresholds that depend on the covariates \citep{northropjonathan2011}.
Such thresholds could retain more observations in regions with lighter conditional tails, where a common threshold may limit power against departures from the null model.
This extension would require controlling the Pareto approximation
bias after projection when the threshold varies across covariate values.

\bibliographystyle{abbrvnat}
\bibliography{references}

\end{document}